\documentclass[final,1p,times]{elsarticle}

\usepackage{epsfig}

\usepackage{amsmath,amsthm,verbatim,amssymb,amsfonts,amscd, graphicx}
\usepackage{latexsym,graphicx, epsfig, bm, epstopdf, float, color}
\usepackage{tikz,siunitx}
\usepackage{scalerel}
\usetikzlibrary{shapes,arrows,patterns}
\usetikzlibrary{shapes.geometric}
\usepackage{pgfplots}
\usetikzlibrary{shapes.geometric,shapes.symbols}
\usetikzlibrary{patterns}
\usetikzlibrary{tikzmark}

\usepackage{caption}
\usepackage{subcaption}
\usepackage{graphics,mathtools,hyperref,cleveref}
\usepackage{dsfont}
\usepackage{framed}
\usepackage[most]{tcolorbox}
\usepackage{xcolor}

\newcommand{\CG}{\text{CG}}
\newcommand{\DG}{\text{DG}}
\newcommand{\EG}{\text{EG}}

\newcommand{\R}{\mathbb{R}}

\newcommand{\x}{\textbf{x}}

\newcommand{\Tau}{\mathcal{T}}
\newcommand{\vertiii}[1]{{\left\vert\kern-0.25ex\left\vert\kern-0.25ex\left\vert #1 
    \right\vert\kern-0.25ex\right\vert\kern-0.25ex\right\vert}}

\theoremstyle{remark}
\newtheorem{remark}{Remark}

\usepackage[textsize=footnotesize]{todonotes}

\journal{Journal of Computational Physics}

\begin{document}

\begin{frontmatter}



\title{Bound-preserving and entropy stable enriched Galerkin methods for nonlinear hyperbolic equations}


\author[DK]{Dmitri Kuzmin} 
\author[FSU]{Sanghyun Lee} 
\author[FSU]{Yi-Yung Yang} 


\affiliation[DK]{organization={Institute of Applied Mathematics, TU Dortmund University},
            addressline={Vogelpothsweg 87}, 
            city={Dortmund},
            postcode={D-44227}, 
            country={Germany}}

\affiliation[FSU]{organization={Department of Mathematics, Florida State University},
            addressline={1017 Academic Way}, 
            city={Tallahassee},
            postcode={32306-4510}, 
            state={FL},
            country={USA}}

\begin{abstract}
  In this paper, we develop monolithic limiting techniques for enforcing nonlinear stability constraints in enriched Galerkin (EG) discretizations of nonlinear scalar hyperbolic equations. To achieve local mass conservation and gain control over the cell averages, the space of continuous (multi-)linear finite element approximations is enriched with piecewise-constant functions. The resulting spatial semi-discretization has the structure of a variational multiscale method. For linear advection equations, it is inherently stable but generally not bound preserving. To satisfy discrete maximum principles and ensure entropy stability in the nonlinear    
  case, we use limiters adapted to the structure of our locally conservative EG method. The cell averages are constrained using a flux limiter, while the nodal values of the continuous component are constrained using a clip-and-scale limiting strategy for antidiffusive element contributions. The design and analysis of our new algorithms build on recent advances in the fields of convex limiting and algebraic entropy fixes for finite element methods. In addition to proving the claimed properties of the proposed approach, we conduct numerical studies for two-dimensional nonlinear hyperbolic problems. The numerical results demonstrate the ability of our limiters to prevent violations of the imposed constraints, while preserving the optimal order of accuracy in experiments with smooth solutions.
\end{abstract}




\end{frontmatter}



\section{Introduction}
\label{sec:intro}

Hyperbolic problems find extensive applications across various fields, playing a crucial role in understanding phenomena such as aquifer contaminant transport~\cite{Bear}, shallow water dynamics~\cite{aizinger2002discontinuous,khan2014modeling}, petroleum engineering for oil and gas production~\cite{EWING1985421, nla.cat-vn721408,Oil-Recovery-Wheeler}, and environmental engineering~\cite{bear2010modeling, karatzas2017developments}. The evolution equations that model various transport phenomena are typically derived from fundamental physical laws, such as conservation principles. Since exact weak solutions of nonlinear hyperbolic problems can be discontinuous and nonunique (even for smooth initial data), computation of physically admissible numerical solutions presents significant challenges. A well-designed discretization method should be sufficiently dissipative to prevent spurious oscillations, nonphysical states and convergence to wrong weak solutions. On the other hand, the levels of numerical viscosity should be kept as low as possible to achieve high resolution.

The framework of algebraic flux correction (AFC) \cite{barrenechea2016,kuzmin2012a,kuzmin2012b} makes it possible to incorporate appropriate constraints, such as discrete conservation/maximum principles and entropy inequalities, into finite element discretizations of scalar conservation laws and hyperbolic systems. Splitting a high-order discretization into a property-preserving low-order part and an antidiffusive remainder, AFC schemes constrain the latter in an adaptive and conservative manner. Examples of such nonlinear stabilization techniques include various generalizations \cite{kuzmin2010a,lohner1987,lohmann2016} of Zalesak's flux-corrected transport (FCT) algorithm \cite{Zalesak-limiter}. The element-based FCT schemes introduced in \cite{cotter2016,lohmann2017} in the context of linear advection problems use convex decompositions to keep the nodal values of finite element approximations in the admissible range under local constraints that are easy to enforce using limiters. For nonlinear hyperbolic systems, an edge-based FCT algorithm of this kind was proposed by Guermond et al. \cite{guermond2018}, who used it to enforce preservation of invariant domains, i.e., of convex admissible sets containing all possible states of exact solutions. The monolithic convex limiting (MCL) strategy introduced in \cite{kuzmin2020} makes it possible to constrain spatial semi-discretizations instead of fully discrete schemes. In contrast to FCT-like predictor-corrector algorithms,  MCL discretizations have well-defined residuals and steady-state solutions. Moreover, they can be equipped with limiter-based entropy fixes \cite{KuHaRu2021,kuzmin2020f} based on Tadmor's entropy stability theory \cite{tadmor2003entropy}. A comprehensive review of the state of the art in the design and analysis of limiting techniques for finite elements can be found in \cite{kuzmin2023}.

In the literature on AFC schemes, discretization in space is typically performed using a continuous Galerkin (CG) or discontinuous Galerkin (DG) method. The advantages of DG approximations include the local conservation property and optimal convergence behavior for linear hyperbolic problems with smooth solutions. However, the number of degrees of freedom is considerably greater than that for a CG approximation using finite elements of the same polynomial degree. As an attractive alternative, we consider an enriched Galerkin (EG) method\footnote{In the original publication \cite{becker2003}, the inventors of EG called it a reduced  $\mathbb{P}_1$-discontinuous Galerkin method.}  \cite{Sun2009,LEE2016,LEE201719}, in which a discontinuous $\mathbb{P}_0$ component $\delta u_h$ is added to the continuous $\mathbb{P}_1$ or $\mathbb{Q}_1$ component $u_h$ of the numerical solution $u_h^{\rm EG}=u_h+\delta u_h$.  The EG formulation can be interpreted as a variational multiscale method in which the evolution equation for $\delta u_h$ serves as a subgrid scale model~\cite{juanes2005}. The favorable properties of the DG-$\mathbb{P}_1$ method are achieved at a lower cost, because the EG version uses less degrees of freedom and produces linear systems that can be solved efficiently~\cite{LEE2016}.
Moreover, the local mass conservation property of EG has been exploited in various applications including two-phase flow~\cite{lee2018enriched}, poroelasticity~\cite{choo2018enriched,lee2023locking}, and thermo-poroelasticity~\cite{yi2024physics}. 

The first bound-preserving EG method for linear advection problems was designed in \cite{Kuzmin-EG} using tailor-made limiters of FCT and MCL type for fluxes and element contributions. In the present paper, we extend the MCL version to nonlinear hyperbolic conservation laws in a manner that enables us to enforce entropy inequalities in addition to maximum principles. For the cell averages that define the DG component $\delta u_h$ of $u_h^{\rm EG}$, we use the flux limiting strategy proposed in \cite{kuzmin2021}. Element contributions to the evolution equation for the CG component $u_h$ are constrained using a clip-and-scale limiting algorithm, which we equip with a built-in entropy fix. Extensive numerical experiments demonstrate the ability of our method to produce physically consistent results while achieving optimal convergence rates. To the best of our knowledge, this is the first successful application of AFC tools to EG discretizations of nonlinear hyperbolic problems.

This paper is organized as follows. In Section~\ref{sec.governing}, we formulate the continuous initial value problem and define the entropy solution, which represents a unique vanishing viscosity limit. In Section~\ref{sec.numerics}, we introduce a baseline EG discretization that features a group finite element (GFE) approximation to the nonlinear flux function. In Section~\ref{sec:algrbraic}, we present the bound-preserving and entropy stable (but low-order accurate) local Lax--Friedrichs-type schemes for the cell averages and CG components of our finite element approximation. In Section~\ref{sec:MCL}, we define the limited anti\-diffusive terms that recover the high-order EG target in smooth regions, while preserving local bounds and maintaining entropy stability of the semi-discrete AFC scheme. In particular, we present and analyze our new entropic clip-and-scale limiter here. In Section~\ref{sec:SSp-RK}, we discretize in time using an explicit strong stability preserving Runge--Kutta (SSP2/Heun) method. In Section~\ref{sec:4:Numerical Results}, we conduct in-depth numerical studies for scalar nonlinear hyperbolic problems. For visualization purposes, we perform constrained $L^2$ projections of EG solutions into the CG space (as in \cite{Kuzmin-EG}). Finally, the results are discussed and conclusions are drawn in Section \ref{sec:conclusions}.

\section{Continuous problem}
\label{sec.governing}

Let $u(\mathbf x, t)\in\mathbb R$ denote the value of a scalar conserved quantity $u$ at a space location $\mathbf x\in\mathbb R^d$, $d \in \{1,2,3\}$ and time $t\geq 0$.  We consider the initial value problem
\begin{subequations}\label{eq:transport}
\begin{alignat}{4}
\dfrac{\partial u}{\partial t} + \nabla \cdot  \mathbf f(u)  &= 0
&&\text{in } \Omega\times \mathbb R_+,\label{eq:transport-a}\\
u(\cdot ,0) &= u_0 \quad &&\text{in }\Omega,
\end{alignat}
\end{subequations}
where $\Omega \subset \mathbb R^d$ is a Lipschitz domain, $\mathbf f(u) := \left(f_1(u),...,f_d(u)\right)$ is the flux function of the (possibly nonlinear) hyperbolic conservation law, and $u_0:\Omega\to\mathbb R$ denotes the initial datum. 

If periodic boundary conditions are imposed on $\partial\Omega$, the formulation \eqref{eq:transport} of the continuous problem is complete.
In the non-periodic case, we prescribe the inflow boundary condition
\begin{equation}
    u = u_\text{in},\quad \text{on }\Gamma^\text{in}\times \mathbb R_+ 
\end{equation}
at the inlet $\Gamma^\text{in} := \{\mathbf x\in\partial \Omega : \mathbf f'(u) \cdot \mathbf n < 0\}$, where $\mathbf f'(u) = \left(f_1'(u),...,f'_d(u)\right)$ is the flux Jacobian and $\mathbf n$ is the unit outward normal to $\partial\Omega$. No boundary condition is imposed on $\Gamma^\text{out} := \partial \Omega \backslash \Gamma^\text{in}$.

Because of hyperbolicity, there exists a convex entropy $\eta(u)$ and an associated entropy flux $\mathbf q(u)$, i.e., a function $\ \mathbf{q}: \mathbb R \to \mathbb R^{ d}$ of $\eta: \mathbb R \to \mathbb R$ such that $\mathbf q'(u) = \eta'(u)\mathbf f'(u)$. If problem (\ref{eq:transport}) has a smooth classical solution $u$, the entropy conservation law
\begin{equation}
    \dfrac{\partial \eta(u)}{\partial t} + \nabla \cdot \mathbf q(u) = 0\qquad \text{in }\Omega \times \mathbb R_+
\end{equation}
can be derived from (\ref{eq:transport-a}) using multiplication by the entropy variable $v(u) := \eta'(u)$, the chain rule, and the above definition of an entropy pair $\left\{\eta(u),\mathbf q(u)\right\}$. In general, the unique vanishing viscosity solution of (\ref{eq:transport}) satisfies a weak form of the entropy inequality \cite{tadmor2003entropy,evans10}
\begin{equation}\label{eq:entropy inequality}
    \dfrac{\partial \eta(u)}{\partial t} + \nabla \cdot \mathbf q(u) \leq 0\qquad \text{in }\Omega \times \mathbb R_+
\end{equation}
for any entropy pair. Hence, entropy is conserved in smooth regions and dissipated at shocks. 



\section{High-order space discretization}
\label{sec.numerics}

\begin{figure}[!h]
\centering
\begin{subfigure}[b]{0.32\textwidth}
\centering
\hspace{-1.6cm}\begin{tikzpicture}
\draw[] (0,0) rectangle (2,2);
\draw (0,1) -- (2,1); 
\draw (0,0.5) -- (2,0.5); 
\draw (0,1.5) -- (2,1.5); 
\draw (2,2) -- (0,2); 
\draw (1,2) -- (1,0); 

\draw (0.5,2) -- (0.5,0); 
\draw (1.5,2) -- (1.5,0); 

\node[] at (0.75,0.75) {$K_e$};

\draw[solid, fill= black] (1,1) circle (0.05);
\draw[solid, fill= black] (0.5,1.) circle (0.05);
\draw[solid, fill= black] (0.5,0.5) circle (0.05);
\draw[solid, fill= black] (1.,0.5) circle (0.05);

\node (a) [] at (1,1) {};
\node (b) [] at (0.5,1) {};
\node (c) [] at (0.5,0.5) {};
\node (d) [] at (1,0.5) {};

\node (N) at (-1,1.) {$\bigcup_{i\in\mathcal N^e }\mathbf x_i$};

\draw (a) edge[->,bend right=25] node {} (N);
\draw (b) edge[->,bend right=15] node {} (N);
\draw (c) edge[->,bend left=15] node {} (N);
\draw (d) edge[->,bend left=5] node {} (N);





\end{tikzpicture}

\caption{$\mathcal N^e$}
\end{subfigure}
\begin{subfigure}[b]{0.32\textwidth}
\centering
\hspace{-1.5cm}\begin{tikzpicture}


\draw [fill=orange,opacity=0.4] (0.5,0.5) rectangle (1.5,1.5);

\node (a) [] at (0.9,0.75) {};
\node (b) [] at (0.9,1.25) {};
\node (c) [] at (1.4,0.75) {};
\node (d) [] at (1.4,1.25) {};

\node (N) at (-1,1.) {$\bigcup_{e\in\mathcal E_i}K_e$};

\draw (a) edge[->] node {} (N);
\draw (b) edge[->] node {} (N);
\draw (c) edge[->,bend left=25] node {} (N);
\draw (d) edge[->,bend right=25] node {} (N);

\draw[] (0,0) rectangle (2,2);
\draw (0,1) -- (2,1); 
\draw (0,0.5) -- (2,0.5); 
\draw (0,1.5) -- (2,1.5); 
\draw (2,2) -- (0,2); 
\draw (1,2) -- (1,0); 

\draw (0.5,2) -- (0.5,0); 
\draw (1.5,2) -- (1.5,0); 

\node[] at (1.15,1.15) {$\mathbf x_i$};

\draw[solid, fill=black] (1,1) circle (0.05);
\draw[solid] (1,0) circle (0.05);
\draw[solid] (0,1) circle (0.05);
\draw[solid] (2,1) circle (0.05);
\draw[solid] (1,2) circle (0.05);
\draw[solid] (0,0) circle (0.05);
\draw[solid] (2,0) circle (0.05);
\draw[solid] (2,2) circle (0.05);
\draw[solid] (0,2) circle (0.05);

\draw[solid] (0.5,0.) circle (0.05);
\draw[solid] (0.5,0.5) circle (0.05);
\draw[solid] (0.5,1.) circle (0.05);
\draw[solid] (0.5,1.5) circle (0.05);
\draw[solid] (0.5,2.) circle (0.05);
\draw[solid] (1.5,0.) circle (0.05);
\draw[solid] (1.5,0.5) circle (0.05);
\draw[solid] (1.5,1.) circle (0.05);
\draw[solid] (1.5,1.5) circle (0.05);
\draw[solid] (1.5,2.) circle (0.05);

\draw[solid] (0,0.5) circle (0.05);
\draw[solid] (0,1.5) circle (0.05);
\draw[solid] (1.,0.5) circle (0.05);
\draw[solid] (1.,1.5) circle (0.05);
\draw[solid] (2.,0.5) circle (0.05);
\draw[solid] (2.,1.5) circle (0.05);
\end{tikzpicture}
\caption{$\mathcal E_i$}
\end{subfigure}
\begin{subfigure}[b]{0.32\textwidth}
\centering
\hspace{-1.8cm}\begin{tikzpicture}


\draw[dashed] (0,0) rectangle (2,2);

\draw [fill=orange,opacity=0.4] (0.5,0.5) rectangle (1.,1.);
\draw [thick, color=red, line width=0.5mm] (0.5,0.5) rectangle (1,1); 

\node (a) [] at (0.75,0.5) {};
\node (b) [] at (0.75,1.) {};
\node (c) [] at (0.6,0.85) {};
\node (d) [] at (1.1,0.75) {};

\node (N) at (-1.2,1.) {$S_{ee'}$};

\draw (d) edge[->,bend left=25] node {} (N);

\draw [dashed] (0,1) -- (2,1); 
\draw [dashed] (0,0.5) -- (2,0.5); 
\draw [dashed] (0,1.5) -- (2,1.5); 
\draw [dashed] (1,2) -- (1,0); 

\draw [dashed] (0.5,2) -- (0.5,0); 
\draw [dashed] (1.5,2) -- (1.5,0); 

\node[] at (0.75,0.75) {$K_e$};

\node[] at (1.25,0.75) {$K_{e'}$};
\node[] at (0.25,0.75) {};

\node[] at (0.75,1.25) {};


\node[] at (0.75,0.25) {};

\end{tikzpicture}
\caption{$S_{ee'}$}
\end{subfigure}

\caption{Visualization of the notation for vertices, cells, and faces of a uniform quadrilateral mesh.}
\label{fig:notation}
\end{figure}
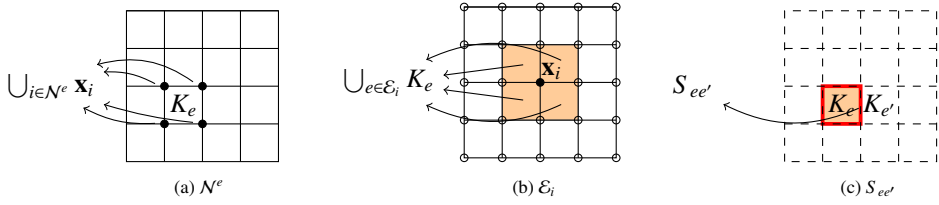

Let $\Tau_h = \{K_e\}_{e=1}^{E_h}$ denote a non-degenerate partition of the domain $\Omega$ into $E_h$ rectangular cells of maximum diameter $h=\max_{1\leq e\leq E_h} h_{K_e}$, where $h_{K_e}$ is the diameter of $K_e$. The vertices of $\Tau_h$ are denoted by $\mathbf x_1,\ldots,\mathbf x_{N_h}$. We store the indices of vertices belonging to a given cell $K_e$ in the integer set $\mathcal N^e$ and the indices of elements that contain a given vertex $\mathbf x_i$ in the integer set $\mathcal E_i$. The set $\mathcal N_i=\bigcup_{e\in\mathcal E_i}\mathcal N^e$ contains the indices of all vertices belonging to at least one cell that contains $\mathbf x_i$. Note that $i\in\mathcal N_i$. The boundary of a cell $K_e$ consists of faces $S_{ee'}$ on which the outward normal $\mathbf n_{ee'}$ is constant. If $S_{ee'} = \partial K_e \cap \partial K_{e'}$ is an internal face, then $e'\in \{1,\ldots,E_h\}\backslash\{e\}$ is the index of an adjacent cell $K_{e'}$. Boundary faces $S_{ee'}\subset\partial \Omega$ are associated with $\partial E_h$ ghost cells and labeled using cell indices $e'\in \{E_h+1,\ldots, \Bar{E}_h\}$, where $\bar E_h:= E_h +\partial E_h$. We store the indices of faces of a cell $K_e$ in the set $\mathcal Z_e$ such that $\partial K_e = \cup_{e'\in\mathcal Z_e} S_{ee'}$. The unit outward normal to $\partial K_e$ is denoted by $\mathbf n_{e}$ if there is no ambiguity. See Figure \ref{fig:notation} for a visual explanation of the above notation.

\subsection{Enriched Galerkin method}
Let $\mathbb Q_{1}(\hat K)$ denote the space of multilinear polynomials $\hat v:\hat K\to\R$ defined on the reference element $\hat K=[0,1]^d$. Using a multilinear mapping $F_e:\hat K\to K_e\in\mathcal T_h$, we construct the space $\mathbb Q_{1}(K_e)$ of polynomials $v:K_e\to\R$ such that $v=\hat v\circ F_e^{-1}$ for some $\hat v\in\mathbb Q_{1}(\hat K)$. 

\begin{figure}[!h]
\centering
\begin{subfigure}[]{0.25\textwidth}
\centering
\begin{tikzpicture}
\draw[thick] (0,0) rectangle (2,2);
\draw (0,1) -- (2,1); 
\draw (2,2) -- (0,2); 
\draw (1,2) -- (1,0); 

\draw[solid,fill=red] (1,1) circle (0.05);
\draw[solid,fill=red] (1,0) circle (0.05);
\draw[solid,fill=red] (0,1) circle (0.05);
\draw[solid,fill=red] (2,1) circle (0.05);
\draw[solid,fill=red] (1,2) circle (0.05);
\draw[solid,fill=red] (0,0) circle (0.05);
\draw[solid,fill=red] (2,0) circle (0.05);
\draw[solid,fill=red] (2,2) circle (0.05);
\draw[solid,fill=red] (0,2) circle (0.05);

\end{tikzpicture}
\caption{CG-$\mathbb{Q}_1$}
\end{subfigure}
\begin{subfigure}[]{0.25\textwidth}
\centering
\begin{tikzpicture}
\draw[thick] (0,0) rectangle (2,2);
\draw (0,1) -- (2,1); 
\draw (2,2) -- (0,2); 
\draw (1,2) -- (1,0); 


\node[mark size=2pt,color=blue] at (0.5,0.5) {\pgfuseplotmark{triangle*}};
\node[mark size=2pt,color=blue] at (0.5+1,0.5+1) {\pgfuseplotmark{triangle*}};
\node[mark size=2pt,color=blue] at (0.5,0.5+1) {\pgfuseplotmark{triangle*}};
\node[mark size=2pt,color=blue] at (0.5+1,0.5) {\pgfuseplotmark{triangle*}};

\end{tikzpicture}
\caption{DG-$\mathbb{Q}_0$}
\end{subfigure}
\begin{subfigure}[]{0.25\textwidth}
\centering
\begin{tikzpicture}

\draw[thick] (0,0) rectangle (2,2);
\draw (0,1) -- (2,1); 
\draw (2,2) -- (0,2); 
\draw (1,2) -- (1,0);

\draw[solid,fill=red] (1,1) circle (0.05);
\draw[solid,fill=red] (1,0) circle (0.05);
\draw[solid,fill=red] (0,1) circle (0.05);
\draw[solid,fill=red] (2,1) circle (0.05);
\draw[solid,fill=red] (1,2) circle (0.05);
\draw[solid,fill=red] (0,0) circle (0.05);
\draw[solid,fill=red] (2,0) circle (0.05);
\draw[solid,fill=red] (2,2) circle (0.05);
\draw[solid,fill=red] (0,2) circle (0.05);

\node[mark size=2pt,color=blue] at (0.5,0.5) {\pgfuseplotmark{triangle*}};
\node[mark size=2pt,color=blue] at (0.5+1,0.5+1) {\pgfuseplotmark{triangle*}};
\node[mark size=2pt,color=blue] at (0.5,0.5+1) {\pgfuseplotmark{triangle*}};
\node[mark size=2pt,color=blue] at (0.5+1,0.5) {\pgfuseplotmark{triangle*}};

\end{tikzpicture}
\caption{EG-$\mathbb{Q}_1$}
\end{subfigure}
\caption{Degrees of freedom for finite element approximations
  on a patch of four cells.}
\label{fig:dof}
\end{figure}
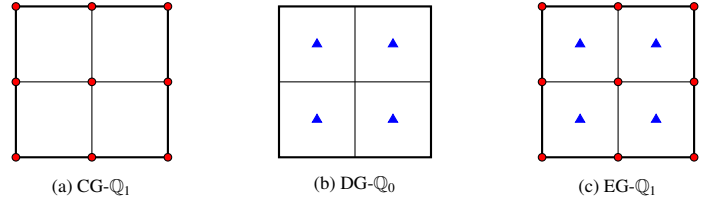

The finite element space CG-$\mathbb Q_1$ of the classical continuous Galerkin method using $\mathbb Q_1$ Lagrange elements on the partition $\Tau_h$ is defined as
\[
V^{\text{CG}}_h := \{v \in L^2(\Omega) : v|_{K_e} \in \mathbb Q_1(K_e)\ \forall K_e \in \Tau_h \} \cap \mathbb C(\bar \Omega),
\]
where $\mathbb C(\bar\Omega)$ denotes the space of functions that are
continuous on $\bar \Omega$. The
DG-$\mathbb Q_0$ space

\begin{equation}
  V^{\DG}_{h} := \{v\in L^2(\Omega)\,:\,
  v|_{K_e}\in\mathbb Q_0(K_e)\ \forall K_e\in \Tau_h\}
\end{equation}
consists of functions that are constant on elements of the partition $\Tau_h$. As already mentioned in the introduction, we approximate exact weak solutions of \eqref{eq:transport} by
\begin{equation}\label{uEG-def}
  u_h^{\EG}=u_h+\delta u_h\in V^{\text{CG}}_h\oplus V^{\DG}_{h},
\end{equation}
where $u_h\in V^{\text{CG}}_h$ and $\delta u_h\in V^{\DG}_{h}$. The
 space of such enriched Galerkin (EG-$\mathbb{Q}_1$)
 approximations (cf. \cite{becker2003,Kuzmin-EG,LEE2016,LEE201719})
is denoted by $V^{\EG}_{h}$. Figure \ref{fig:dof}
shows the local degrees of freedom (DOFs) for CG-$\mathbb Q_1$,
DG-$\mathbb Q_0$, and EG-$\mathbb{Q}_1$ approximations on a
patch consisting of four square cells.

\begin{remark}
  We use quadrilateral meshes in this work but piecewise-linear (CG-$\mathbb{P}_1$) approximations $u_h$ on triangles
  can be enriched by piecewise-constant functions $\delta u_h\in V^{\DG}_{h}$ similarly. The linear stability analysis performed in \cite{becker2003} for EG discretizations of convection-diffusion equations is valid only for triangular meshes. However, the piecewise-constant enrichment of CG-$\mathbb Q_1$ approximations also provides the local conservation property and has a stabilizing effect.
\end{remark}

The space $V^{\CG}_{h}$ is spanned by $N_h$ global basis functions $\{\varphi_j\}_{j=1}^{N_h}$ such that $\varphi_j\in V^{\CG}_{h}$ and $\varphi_j(\mathbf{x}_i)=\delta_{ij}$ for all $i,j\in\{1,\ldots,N_h\}$. Hence, the CG component can be written as
 \begin{equation}
 u_h= \sum_{j=1}^{N_h}u_j\varphi_j\in V^\text{CG}_{h},
 \label{eqn:u_i}
 \end{equation}
where $u_j=u_h(\mathbf{x}_j)$ is the degree of freedom associated with $\varphi_j$. Similarly, the DG component 
\begin{equation}\label{deltauh}
\delta u_h = \sum_{e=1}^{E_h}\delta u_e \chi_{e}\in V^\text{DG}_{h}
\end{equation}
is defined by the coefficients $\delta u_e$ of $E_h$ characteristic functions $\chi_e$ such that $\chi_{e} = 1$ in $K_e$ and $\chi_{e} = 0$ otherwise.  The piecewise-constant basis functions $\chi_e$ span the space $V^\text{DG}_{h}$.
\medskip
  
Using the EG method to discretize the nonlinear problem \eqref{eq:transport}, we seek a numerical solution $u^\text{EG}_h \in V^\text{EG}_{h}$ that admits the decomposition \eqref{uEG-def} and satisfies
the DG weak form \cite{kuzmin2023}
\begin{align}\nonumber
    \sum_{e=1}^{E_h}\int_{K_e} w^\text{EG}_h \dfrac{\partial u^\text{EG}_h}{\partial t} \mathrm d\mathbf x 
    &+ \sum_{e=1}^{E_h}\int_{\partial K_e} w^\text{EG}_h [\mathcal F(u_{h,e}^-,u_{h,e}^+;\mathbf n_{e}) - \mathbf f(u_{h,e}^\text{EG})\cdot \mathbf n_e] \mathrm ds\\ 
    &+ \sum_{e=1}^{E_h}\int_{K_e} w^\text{EG}_h\nabla \cdot \mathbf f(u_h^\text{EG}) \mathrm d\mathbf x = 0\label{eq:consistentEG}
\end{align}
of (\ref{eq:transport-a})
for all test functions $ w^{\EG}_h \in V^\text{EG}_{h}$. We define the
local Lax-Friedrichs  (LLF) flux
\begin{equation}\label{eq:LLF}
    \mathcal F(u_L, u_R;\mathbf n) := \dfrac{\mathbf f(u_R)+\mathbf f(u_L)}{2}\cdot \mathbf n - \dfrac{\lambda_{LR}}{2}(u_R - u_L)
\end{equation}
using the maximum wave speed \cite{kuzmin2021}
\begin{equation}\label{eq:wave-speed}
    \lambda_{LR} =\lambda_{\max}(u_L,u_R;\mathbf n) := \max_{\omega \in[0,1]} |\mathbf {f}'(\omega u_R + (1-\omega)u_L)\cdot \mathbf n|.
\end{equation}
  The internal and external states of the local Riemann problem are defined by
\begin{equation*}
  u_{h,e}^-=u_h+\delta u_e,\qquad
   u_{h,e}^+= 
    \begin{cases}
      \  u_\text{in}\qquad &\text{on }\Gamma^\text{in},\\
       \  u_h+\delta u_e\qquad &\text{on }\Gamma^\text{out},\\
      \   u_h+\delta u_{e'}\qquad &\text{on }S_{ee'} = \partial K_e \cap \partial K_{e'},\ e'\le E_h.
    \end{cases}
\end{equation*}

Integration by parts for the last term on the left-hand side of (\ref{eq:consistentEG}) yields 
\begin{equation}\label{eq:consisten formula 2}
       \sum_{e=1}^{E_h}\int_{K_e} w^\text{EG}_h \dfrac{\partial u^\text{EG}_h}{\partial t} \mathrm d\mathbf x 
       + \sum_{e=1}^{E_h}\int_{\partial K_e} w^\text{EG}_h \mathcal F(u_{h,e}^-,u_{h,e}^+;\mathbf n_e) \mathrm ds - \sum_{e=1}^{E_h}\int_{K_e} \nabla  w^\text{EG}_h\cdot \mathbf f(u_h^\text{EG}) \mathrm d\mathbf x = 0.
\end{equation}

Any test function $w^{\EG}_h \in V^\text{EG}_{h}$ can be written as
$w^{\EG}_h=w_h+\delta w_h$, where  $w_h \in V_h^{\rm CG}$ and
$\delta w^{\EG}_h \in V_h^{\rm DG}$. In particular,
$w_h$ and $\delta w_h$ are admissible test
functions. Thus
\begin{align}\nonumber
       \sum_{e=1}^{E_h}\int_{K_e} w_h \dfrac{\partial u^\text{EG}_h}{\partial t} \mathrm d\mathbf x 
       &+ \sum_{e=1}^{E_h}\int_{\partial K_e\cap\partial\Omega} w_h \mathcal F(u_{h,e}^-,u_{h,e}^+;\mathbf n_e) \mathrm ds\\
       &- \sum_{e=1}^{E_h}\int_{K_e} \nabla  w_h\cdot \mathbf f(u_h^\text{EG}) \mathrm d\mathbf x = 0\qquad
       \forall w_h\in V_h^{\CG},\label{subproblemCG}\\
\sum_{e=1}^{E_h}\int_{K_e} \delta w_h \dfrac{\partial u^\text{EG}_h}{\partial t} \mathrm d\mathbf x 
&+ \sum_{e=1}^{E_h}\int_{\partial K_e} \delta  w_h \mathcal F(u_{h,e}^-,u_{h,e}^+;\mathbf n_e) \mathrm ds { = 0}\qquad \forall  \delta  w_h\in V_h^{\DG}.
\label{subproblemDG}
\end{align}
Note that the surface integration in \eqref{subproblemCG} is restricted to boundary faces of $\mathcal T_h$. The integrals over internal faces $S_{ee'}=\partial K_e\cap \partial K_{e'}$ cancel out for continuous test functions $w_h$.

As noticed by Becker et al. \cite{becker2003}, the representation \eqref{uEG-def} of $u^\text{EG}_h \in V^\text{EG}_{h}$ is nonunique because any globally constant function can be represented exactly in $V_h^{\CG}$ and in $V_h^{\DG}$ alike. To ensure the uniqueness of $\delta u_h$ defined by \eqref{deltauh}, we impose the additional constraint \cite{Kuzmin-EG}
\begin{equation}\label{eq:massless}
    \delta u_e = \dfrac{1}{|K_e|}\int_{K_e} (u^\text{EG}_h - u_h)\mathrm d\mathbf x = U_e - \bar u_e,\qquad e = 1,\ldots, E_h,
\end{equation}
where
\begin{equation}\label{eq:averages}
U_e=\bar u_e = \frac{1}{|K_e|}\int_{K_e}u_h^{\EG} \mathrm d\mathbf x,\qquad
\bar u_e = \frac{1}{|K_e|}\int_{K_e}u_h \mathrm d\mathbf x,\qquad
|K_e|=\int_{K_e}1\mathrm d\mathbf x.
\end{equation}

In the remainder of this section, we derive evolution equations for
the  cell averages $U_e$ of $u_h^{\EG}$ and the
nodal values $u_j$ 
of $u_h$ defined by \eqref{eqn:u_i}. The evolution of the DG component
$\delta u_h=u_h^{\EG}-u_h$ is then determined by our convention \eqref{eq:massless}, which implies that
$\delta u_h$ is \emph{massless} \cite{Kuzmin-EG}. 

\begin{remark}
  The analogy with variational multiscale methods 
  for conservation laws (cf.~\cite{juanes2005})
makes it possible to interpret equations \eqref{subproblemDG} and
\eqref{eq:massless} as a subgrid scale model for $\delta u_h$. 
\end{remark}

\subsection{High-order semi-discrete problem for EG cell averages}
    \label{sec:HOcell}

Using  test functions $\delta w_h\in \{\chi_1,\ldots,\chi_{E_h}\}$ in the semi-discrete weak form
\eqref{subproblemDG}, we find that the cell averages
of the EG solution $u_h^\text{EG}$ must satisfy the
local conservation laws
\begin{equation}\label{eq:cell_average}
  |K_e|\dfrac{\mathrm dU_e}{\mathrm dt} =
  -\int_{\partial K_e} \mathcal F(u_{h,e}^-,u_{h,e}^+;\mathbf n_e)\mathrm ds
=:q_e^H,\qquad e=1,\ldots,E_h,
\end{equation}
where
  $$
  u_{h,e}^-=u_h+(U_e-\bar u_e),\qquad u_{h,e'}^+=
   \begin{cases}
      \  u_\text{in}\qquad &\text{on }\Gamma^\text{in},\\
       \ u_h+(U_e-\bar u_e) &\text{on }\Gamma^\text{out},\\
      \ u_h+(U_{e'}-\bar u_{e'}) &\text{on }S_{ee'} = \partial K_e \cap \partial K_{e'},\ e'\le E_h.
    \end{cases}
   $$
   The right-hand side $q^H_e$ of equation \eqref{eq:cell_average}
   can be expressed in terms of the face-averaged fluxes
\begin{equation}\label{HQ1aver}
  H_{ee'}^{\mathbb{Q}_1}= \dfrac{1}{|S_{ee'}|} \int_{S_{ee'}} \mathcal F(u^\text{EG}_{h,e},u^\text{EG}_{h,e'};\mathbf n_{ee'} )\mathrm ds 
  \end{equation}
as follows:
\begin{equation*}
  q^H_e =  -\sum_{e'\in\mathcal Z_e}|S_{ee'}| H^{\mathbb Q_1}_{ee'}.
\end{equation*}
Throughout this paper, the superscript $H$ refers to a
`high-order' approximation. The superscript $L$ is reserved for the
`low-order' approximations to be presented in Section
\ref{sec:algrbraic}.
\smallskip

Let $\mathbf U = (U_e)_{e=1}^{E_h}$ denote the vector of EG cell averages
and $\mathbf q^H= (q^H_e)_{e=1}^{E_h}$ the vector of discrete right-hand sides.
By definition, $\mathbf q^H=\mathbf q^H(\mathbf U,\mathbf u)$
depends not only on $\mathbf U$ but also on the  vector
$\mathbf u = (u_i)_{i=1}^{N_h}$ of CG nodal values. Introducing
the finite volume mass matrix $\Bar{M} = \text{diag}(|K_e|)_{e=1}^{E_h}$,
we write the system of ordinary differential equations
 \eqref{eq:cell_average} in the matrix form
\begin{equation}\label{eq:HOsysDG}
    \Bar{M} \dfrac{\mathrm d\mathbf U}{\mathrm dt} = \mathbf q^H(\mathbf U, \mathbf u),
\end{equation}
which will be referred to as the high-order
(HO) semi-discrete problem for EG cell averages.

\subsubsection{High-order semi-discrete problem for CG nodal values}

Using the test functions $w_h^{\rm EG}\in \{\varphi_1,\ldots,\varphi_{N_h}\}$
in \eqref{eq:consistentEG}, we obtain the semi-discrete equations
\begin{align}\notag
    \sum_{e\in\mathcal E}\int_{K_e}\varphi_i \dfrac{\partial u^\text{EG}_h}{\partial t} \mathrm d\mathbf x 
    &+\sum_{e\in\mathcal E_i}\int_{\partial K_e\cap\partial\Omega}\varphi_i
        \big(\mathcal F(u_{h,e}^-,u_{h,e}^+;\mathbf n_{e})-\mathbf f(u_h^{\EG}\big)\cdot
          \mathbf n_e)\mathrm d s\\ 
    &+\sum_{e\in\mathcal E_i}\int_{K_e}\varphi_i\nabla\cdot \mathbf f(u_h^\text{EG}) \mathrm d\mathbf x = 0,\qquad i = 1,\ldots, N_h
\label{eqn:before_taylor}       
\end{align}
for the nodal values of the CG component $u_h$. The union $\bar\Omega_i=\cup_{e\in\mathcal E_i} K_e$ of elements that contain the vertex $\mathbf x_i$ represents the compact support of the basis function $\varphi_i$.

To derive a high-order semi-discrete scheme that is better suited for limiting, we first modify the last term on the left-hand side of \eqref{eqn:before_taylor} using the linear Taylor approximation (cf.~\cite{juanes2005})  
\begin{equation}
    \mathbf f(u_h^\text{EG})=\mathbf f(u_h+\delta u_h)  \approx \mathbf f(u_h) + \mathbf f'(u_h)\delta u_h.
\label{eqn:taylor}    
\end{equation}
Next, we replace $\mathbf f(u_h)$ by the \emph{group finite element} (GFE) interpolant
\cite{kuzmin2012b,FLETCHER1983225,selmin1996}
\begin{equation}\label{eq:GFE}
  \mathbf f_h(u_h)= \sum_{j=1}^{N_h}\mathbf f(u_j)\varphi_j.
\end{equation}
This modification corresponds to using inexact nodal quadrature for
volume integrals. For
$$
\int_{\partial K_e\cap \partial \Omega} \varphi_i\big(\mathcal F(u_{h,e}^-,u_{h,e}^+;\mathbf n_e)-\mathbf f(u_h^{\EG})\cdot
          \mathbf n_e\big)\mathrm d s,
$$
we use the `lumped' GFE approximation  (cf. \cite[Eq. (3.118)]{kuzmin2023})
$$
\sum_{e'\in\mathcal E_e'}\sigma_{i,ee'}
[\mathcal F(u_i,\hat u_{i,ee'};\mathbf n_{ee'})-\mathbf f(u_i)\cdot\mathbf n_{ee'}],
$$
where $\mathcal E_e'=\{e'\in\mathcal Z_e\,:\, e'>E_h\}$ is the
set of second subscripts of $S_{ee'}\subset\partial K_e\cap\partial\Omega$ and
$$
  \sigma_{i,ee'}=\int_{S_{ee'}}\varphi_i\mathrm d s,\qquad
\hat u_{i,ee'}=\frac{1}{\sigma_{i,ee'}}\int_{S_{ee'}}\varphi_iu_{h,e}^+\mathrm d s.
$$

Finally, we modify the time derivative term by using the convenient approximation 
$$
\int_{K_e}\varphi_i \dfrac{\partial u^\text{EG}_h}{\partial t} \mathrm d\mathbf x
\approx \left(\int_{K_e}\varphi_i
\mathrm d\mathbf x\right)\frac{\mathrm d u_i}{\mathrm d t}+
\int_{K_e}\varphi_i(\dot u_h-\dot u_i)\mathrm d\mathbf x,
$$
where $u_i$ is the time-dependent value of the
CG component $u_h$ at the vertex $\mathbf x_i$.
Following Kuzmin et al. \cite{Kuzmin-EG}, 
we reconstruct the nodal values of 
$\dot u_h=\sum_{j=1}^{N_h}\dot u_j\varphi_j$ as follows:
\begin{equation}
\dot u_i=\frac{1}{|\Omega_i|}\sum_{e\in\mathcal E_i}|K_e|
\frac{\mathrm d U_e}{\mathrm d t}
=\frac{1}{|\Omega_i|}\sum_{e\in\mathcal E_i}q_e^H,\qquad
|\Omega_i|=\sum_{e\in\mathcal E_i}|K_e|.\label{eq:timeder}
\end{equation}
The resulting system of evolution equations for the CG nodal values reads
\begin{align}\notag
    \left(\sum_{e\in\mathcal E}\int_{K_e}\varphi_i
    \mathrm d\mathbf x\right)\frac{\mathrm d u_i}{\mathrm d t}
    &=\sum_{e\in\mathcal E}\int_{K_e}\varphi_i(\dot u_i-\dot u_h)\mathrm d\mathbf x
    -\sum_{e\in\mathcal E_i}\sum_{e'\in\mathcal E_e'}\sigma_{i,ee'}
    [\mathcal F(u_i,\hat u_{i,ee'};\mathbf n_{ee'})-\mathbf f(u_i)\cdot\mathbf n_{ee'}]
   \notag \\
    &- \sum_{e\in\mathcal E_i}\sum_{j\in\mathcal N^e}\left(
    \int_{K_e}\varphi_i\nabla\varphi_j\mathrm d\mathbf x\right)
    \cdot\mathbf{f}(u_j)\notag\\
    &+\sum_{e\in\mathcal E_i}  \delta u_e
    \int_{K_e}\nabla \varphi_i\cdot\mathbf f'(u_h)\mathrm d\mathbf x,\qquad i = 1,\ldots, N_h
\label{eqn:lumpedHOscheme}       
\end{align}
and represents a second-order perturbation of \eqref{eqn:before_taylor}.
While it is possible to correct the corresponding errors in the limiting step,
  doing so would make the algorithms to be presented in Section \ref{sec:MCL}
  more complicated without having any significant positive impact on the
  accuracy of numerical solutions.  Therefore, we use 
  \eqref{eqn:lumpedHOscheme} rather than
  \eqref{eqn:before_taylor} as the high-order target for limiting.
  \smallskip
  
  To write system \eqref{eqn:lumpedHOscheme} in a matrix form, we need to
  introduce some further notation. Let
 \begin{align*}
     m_{ij} &= \sum_{e\in\mathcal E_i\cap \mathcal E_j} m^e_{ij},\qquad m^e_{ij} = \int_{K_e} \varphi_i\varphi_j \mathrm d\mathbf x,\\
     m_i &= \sum_{e\in\mathcal E_i} m^e_i, \qquad m^e_i  = \sum_{j\in\mathcal N^e} m^e_{ij} =\int_{K_e}\varphi_i \mathrm d\mathbf x,\\
     \mathbf c_{ij} &= \sum_{e\in\mathcal E_i\cap \mathcal E_j} \mathbf c^e_{ij},\qquad
    \mathbf c^e_{ij} = \int_{K_e}\varphi_i\nabla \varphi_j\mathrm d\mathbf x
\end{align*}
  denote the coefficients of the consistent mass matrix
  $M_C = (m_{ij})_{i,j=1}^{N_h}$, of its lumped counterpart
  $M_L = \text{diag}(m_i)_{i=1}^{N_h}$, and of the discrete gradient
  operator  $\mathbf C = (\mathbf c_{ij})_{i,j=1}^{N_h}$. Note that
  $$
  \mathbf c_{ji}=-\mathbf c_{ij}+\sum_{e\in\mathcal E_i}
  \int_{\partial K_e\cap\partial\Omega}\varphi_i\varphi_j \mathbf{n}\mathrm ds,\qquad
  \mathbf c_{ji}^e=-\mathbf c_{ij}^e+\int_{\partial K_e}\varphi_i\varphi_j
 \mathbf{n}\mathrm ds.
 $$
Adopting the above notation, we  write the evolution equations
 \eqref{eqn:lumpedHOscheme} in the equivalent form
\begin{align}\notag
m_i\frac{\mathrm d u_i}{\mathrm d t}
&=\sum_{e\in\mathcal E_i}
\sum_{j\in\mathcal N^e}m_{ij}^e(\dot u_i-\dot u_j)
    -\sum_{e\in\mathcal E_i}\sum_{e'\in\mathcal E_e'}\sigma_{i,ee'}
    [\mathcal F(u_i,\hat u_{i,ee'};\mathbf n_{ee'})-\mathbf f(u_i)\cdot\mathbf n_{ee'}]
   \notag \\
    &- \sum_{e\in\mathcal E_i}\sum_{j\in\mathcal N^e}\mathbf c_{ij}^e
    \cdot\mathbf{f}(u_j)
    +\sum_{e\in\mathcal E_i}  (U_e-\bar u_e)
    \int_{K_e}\nabla \varphi_i\cdot\mathbf f'(u_h)\mathrm d\mathbf x=:g_i^H,\qquad i = 1,\ldots, N_h.
\label{eqn:lumpedHOscheme2}       
\end{align}
The matrix form of this high-order semi-discrete problem for CG
nodal values is given by
\begin{equation}\label{eq:HOsysCG}
M_L\frac{\mathrm d\mathbf u}{\mathrm d t}=\mathbf g^H(\mathbf U,\mathbf u),
\end{equation}   
where
$\mathbf g^H=(g_i^H)_{i=1}^{N_h}$ is the global vector containing the 
right-hand sides of system \eqref{eqn:lumpedHOscheme2}.

\section{Algebraic splitting}\label{sec:algrbraic}

In Section~\ref{sec.numerics}, we have shown that the 
cell averages $U_e$ and CG nodal values $u_i$ of a high-order EG
approximation $u_h^{\rm EG}$ to the  solution of problem \eqref{eq:transport}
satisfy the nonlinear system
\begin{equation}\label{eq:EG-HO-both}
    \begin{bmatrix}
       \bar{M} & 0 \\
        0 & M_L 
    \end{bmatrix}\dfrac{\mathrm d}{\mathrm dt}\begin{bmatrix}
         \mathbf U  \\
         \mathbf u 
    \end{bmatrix} = \begin{bmatrix}
    {\mathbf q}^H(\mathbf U,\mathbf u)\\
    {\mathbf g}^H(\mathbf U,\mathbf u)
    \end{bmatrix}
\end{equation}
of coupled semi-discrete subproblems \eqref{eq:HOsysDG}
and \eqref{eq:HOsysCG} for $\mathbf U=(U_e)_{e=1}^{E_h}$ and
 $\mathbf u=(u_i)_{i=1}^{N_h}$.
The DG components $\delta u_e=U_e-\bar u_e$
are uniquely determined by the cell averages of $u_h^{\rm EG}$ and $u_h$.

To prepare the ground for the derivation of 
limiting techniques in Section~\ref{sec:MCL}, we need to split the high-order
system \eqref{eq:EG-HO-both} into a property-preserving low-order
part and an antidiffusive correction term. The general 
form of such an algebraic splitting is given by \cite{Kuzmin-EG}
\begin{equation}\label{eq:EG-HO-equiv}
    \begin{bmatrix}
       \bar{M} & 0 \\
        0  & M_L 
    \end{bmatrix}\dfrac{\mathrm d}{\mathrm dt}\begin{bmatrix}
         \mathbf U  \\
         \mathbf u 
    \end{bmatrix} = \begin{bmatrix}
    {\mathbf q}^L(\mathbf U) + \bar{\mathbf f}(\mathbf U,\mathbf u)\\
    {\mathbf g}^L(\mathbf u) + \mathbf{f}(\mathbf U,\mathbf u)  
    \end{bmatrix}.
\end{equation}
The components of the vectors $\bar{\mathbf f}=(\bar f_e)_{e=1}^{E_h}$ and
$\mathbf f=(f_i)_{i=1}^{N_h}$ should admit decompositions
$$
\bar f_e=\sum_{e'\in\mathcal Z_e}F_{ee'},\qquad f_i=\sum_{e\in\mathcal E_i}f_i^e
$$
into numerical fluxes $F_{ee'}$ and element contributions $f_i^e$ such that
\begin{equation}\label{eq:zerosums}
F_{ee'}+F_{e'e}=0,\qquad \sum_{i\in\mathcal N^e}f_i^e=0.
\end{equation}
These zero sum properties must be preserved by limiters to guarantee
discrete conservation.

The right-hand sides ${\mathbf q}^L=(q_e^L)_{e=1}^{E_h}$ and
$\mathbf g^L=(g_i^L)_{i=1}^{E_h}$ of the
low-order (LO) semi-discrete problems for  $\mathbf U$
and $\mathbf u$ must ensure the validity of all relevant constraints for a
solution of
\begin{equation}\label{eq:EG-LO}
    \begin{bmatrix}
       \bar{M} & 0 \\
        0  & M_L 
    \end{bmatrix}\dfrac{\mathrm d}{\mathrm dt}\begin{bmatrix}
         \mathbf U  \\
         \mathbf u 
    \end{bmatrix} = \begin{bmatrix}
    {\mathbf q}^L(\mathbf U)\\
    {\mathbf g}^L(\mathbf u) 
    \end{bmatrix}.
\end{equation}
In the remainder if this section, we present a splitting that provides the above properties.

\subsection{Low-order semi-discrete problem for EG cell averages}

A natural choice of the low-order scheme for $\mathbf U$
is the finite volume LLF method \cite{kuzmin2023,kuzmin2021}
\begin{equation}\label{eq:cell-centered-fv-low}
  |K_e| \dfrac{\mathrm dU_e}{\mathrm dt} =
  -\sum_{e'\in\mathcal Z_e} \int_{S_{ee'}}\mathcal F(U_e^-,U_{e}^+; \mathbf n_{ee'})\mathrm ds
  =:q_e^L,\qquad e = 1,\ldots, E_h,
\end{equation}
in which the LLF fluxes depend on the internal state $U_e^-=U_e$ and the external state
\begin{equation*}
    U_{e}^+ = \begin{cases}
      u_\text{in}\qquad &\text{on }
      \Gamma^\text{in}\\
        U_e\qquad &\text{on }\Gamma^\text{out}\\
        U_{e'}\qquad &\text{on }S_{ee'} = \partial K_e \cap \partial K_{e'},\ e'\le E_h.
    \end{cases}
\end{equation*}
This approximation is known to be bound-preserving and entropy stable
\cite{kuzmin2023,kuzmin2021}. Therefore, it is widely used as the LO component of
modal DG methods equipped with flux and/or slope limiters.
The matrix form of the corresponding
semi-discrete problem reads
\begin{equation}\label{eq:LOsysDG}
    \Bar{M} \dfrac{\mathrm d\mathbf U}{\mathrm dt} = \mathbf q^L(\mathbf U).
\end{equation}
In contrast to $\mathbf q^H(\mathbf U,\mathbf u)$,
the right-hand side of \eqref{eq:LOsysDG}
is independent of $\mathbf u$. Therefore,
the LO cell averages evolve independently from
the LO nodal values of the CG component $u_h$.

The difference between the right-hand sides of equations \eqref{eq:cell_average} and
\eqref{eq:cell-centered-fv-low} is given by
\begin{align*}
\bar f_e=q_e^H-q_e^L&=
\sum_{e'\in\mathcal Z_e} \int_{S_{ee'}}[\mathcal F(U_e^-,U_e^+; \mathbf n_{ee'})
  -\mathcal F(u_{h,e}^-,u_{h,e}^+;\mathbf n_e)]\mathrm ds\\
&=\sum_{e'\in\mathcal Z_e}|S_{ee'}|(H^{\mathbb Q_0}_{ee'}-H^{\mathbb Q_1}_{ee'})
=\sum_{e'\in\mathcal Z_e}F_{ee'},
\end{align*}
where $H^{\mathbb Q_0}_{ee'}$ is the
low-order counterpart of the LLF flux $H^{\mathbb Q_1}_{ee'}$ defined by
\eqref{HQ1aver} and 
\begin{equation}\label{eq:antidiffDG}
F_{ee'}=|S_{ee'}|(H^{\mathbb Q_0}_{ee'}-H^{\mathbb Q_1}_{ee'})=-F_{e'e}
\end{equation}
are the raw antidiffusive fluxes that constitute the component
$\bar f_e$ of the vector $\bar{\mathbf{f}}(\mathbf U,\mathbf u)$
in \eqref{eq:EG-HO-equiv}.

\subsubsection{Low-order semi-discrete problem for CG nodal values}

An algebraic CG version of the LLF method
\eqref{eq:cell-centered-fv-low} approximates the
HO scheme \eqref{eqn:lumpedHOscheme2} by \cite{kuzmin2020g}
\begin{align}\notag
m_i\frac{\mathrm d u_i}{\mathrm d t}
&= -\sum_{e\in\mathcal E_i}\sum_{e'\in\mathcal E_e'}\sigma_{i,ee'}
    [\mathcal F(u_i,\hat u_{i,ee'};\mathbf n_{ee'})-\mathbf f(u_i)\cdot\mathbf n_{ee'}]
   \notag \\
    &+\sum_{e\in\mathcal E_i}\sum_{j\in\mathcal N^e}[d_{ij}^eu_j-\mathbf c_{ij}^e
    \cdot\mathbf{f}(u_j)]
    =:g_i^L,\qquad i = 1,\ldots, N_h.
\label{eqn:LOschemeCG}       
\end{align}
The LLF graph viscosity coefficients
\begin{equation}\label{eq:d_e_ij}
d^e_{ij} = 
\begin{cases}
      \max (\lambda_{ij}^e |\mathbf c^e_{ij}|,\lambda_{ji}^e |\mathbf c^e_{ji}|)&\text{if }\ j\in\mathcal N^e,\ j\ne i,\\
       -\sum_{k\in \mathcal N^e \backslash\{i\}} d^e_{ik}&\text{if }\ j\in\mathcal N^e,\ j=i,\\
       0&\text{otherwise}
\end{cases}
\end{equation}
depend on $\lambda_{ij}^e
= \lambda_{\max}(u_i, u_j;\mathbf{n}_{ij})$,
where $\mathbf{n}_{ij}=\frac{\mathbf c^e_{ij}}{|\mathbf c^e_{ij}|}$ and
$\lambda_{\max}(u_L,u_R;\mathbf{n})$ is defined as in 
\eqref{eq:wave-speed}.
  The matrix form of the 
semi-discrete scheme \eqref{eqn:LOschemeCG} is given by
\begin{equation}\label{eq:LOsysCG}
   M_L \dfrac{\mathrm d\mathbf u}{\mathrm dt} = \mathbf g^L(\mathbf u).
\end{equation}
Since the right-hand side $\mathbf g^L(\mathbf u)$
is independent of $\mathbf U$, so is a solution $\mathbf{u}$ of this LO system. 

A comparison of the spatial semi-discretizations
\eqref{eqn:lumpedHOscheme2} and \eqref{eqn:LOschemeCG} reveals that
\begin{align*}
  f_i=g_i^H-g_i^L&=\sum_{e\in\mathcal E_i}
  \sum_{j\in\mathcal N^e}[m_{ij}^e(\dot u_i-\dot u_j)+d_{ij}^e(u_i-u_j)]\\&
   +\sum_{e\in\mathcal E_i} (U_e-\bar u_e)
   \int_{K_e}\nabla \varphi_i\cdot\mathbf f'(u_h)\mathrm d\mathbf x=
\sum_{e\in\mathcal E_i}f_i^e.
\end{align*}
The individual element contributions $f_i^e$
to this component of $\mathbf{f}(\mathbf U,\mathbf u)$
are defined by
\begin{equation}\label{eq:antidiffCG}
  f_i^e=\sum_{j\in\mathcal N^e}[m_{ij}^e(\dot u_i-\dot u_j)+d_{ij}^e(u_i-u_j)]
   + (U_e-\bar u_e)
   \int_{K_e}\nabla \varphi_i\cdot\mathbf f'(u_h)\mathrm d\mathbf x.
\end{equation}
To prove the zero sum property required in \eqref{eq:zerosums}, we first
notice that $m_{ij}^e=m_{ji}^e$ and $d_{ij}^e=d_{ji}^e$ by definition. The Lagrange basis
functions $\varphi_i$ of the CG-$\mathbb{Q}_1$ space $V_h$ form a partition of unity, i.e.,
$\sum_{i\in\mathcal N^e}\varphi_i\equiv 1$. It follows that  
$\sum_{i\in\mathcal N^e}\nabla\varphi_i\equiv \mathbf 0$ and
$\sum_{i\in\mathcal N^e}f_i^e=0$, as desired.

\section{Property-preserving limiters}\label{sec:MCL}

Given the algebraic splitting \eqref{eq:EG-HO-equiv} of
the high-order spatial semi-discretization \eqref{eq:EG-HO-both}, we
can proceed to constructing a property-preserving EG
approximation of the form
\begin{equation}\label{eq:EG-MCL}
    \begin{bmatrix}
       \bar{M} & 0 \\
        0  & M_L 
    \end{bmatrix}\dfrac{\mathrm d}{\mathrm dt}\begin{bmatrix}
         \mathbf U  \\
         \mathbf u 
    \end{bmatrix} = \begin{bmatrix}
    {\mathbf q}^L(\mathbf U) + \bar{\mathbf f}^*(\mathbf U,\mathbf u)\\
    {\mathbf g}^L(\mathbf u) + \mathbf{f}^*(\mathbf U,\mathbf u)  
    \end{bmatrix},
\end{equation}
where the vectors $\bar{\mathbf f}^*=(\bar f_e^*)_{e=1}^{E_h}$ and
$\mathbf f^*=(f_i^*)_{i=1}^{N_h}$ are composed from
$$
\bar f_e^*=\sum_{e'\in\mathcal Z_e}F_{ee'}^*,\qquad
f_i^*=\sum_{e\in\mathcal E_i}f_i^{e,*}.
$$
The constrained approximations $F_{ee'}^*\approx F_{ee'}$ and
 $f_i^{e,*}\approx f_i^e$
should satisfy the zero sum conditions
\begin{equation}\label{eq:zerosumsMCL}
F_{ee'}^*+F_{e'e}^*=0,\qquad \sum_{i\in\mathcal N^e}f_i^{e,*}=0.
\end{equation}
Additionally, the limiters to be presented below will ensure the validity
of local discrete maximum principles (DMPs) and entropy stability of the semi-discrete schemes for $\mathbf U$ and $\mathbf u$.

\subsection{Limiting criteria}

A closed interval $\mathcal G=[u^{\min},u^{\max}]$ represents
an invariant domain of the scalar hyperbolic problem
\eqref{eq:transport} if preservation of the global bounds
$u^{\min}$ and $u^{\max}$ can be shown for exact solutions.
We call the spatial semi-discretization \eqref{eq:EG-MCL} invariant
  domain preserving (IDP) if it produces admissible states
$U_e(t)\in\mathcal G$ 
and $u_i(t)\in\mathcal G$ at any time $t>0$. We need
\eqref{eq:EG-MCL} to be IDP and, moreover,
local extremum diminishing (LED)
in the sense that the time derivatives of $U_e(t)$ and $u_i(t)$ are
nonpositive/nonnegative at a local maximum/minimum
\cite{jameson1993,jameson1995,kuzmin2012a,kuzmin2023}.

Suppose that there exist bounded solution-dependent
coefficients $A_e\ge 0$
and $a_i\ge 0$ such that
\begin{align}\label{eq:IDPcritDG}
  \dfrac{\mathrm dU_e}{\mathrm dt}&=
  \frac{q_e^L + \bar{f}_e^*}{|K_e|}= A_e(\bar U_e^*-U_e),\qquad e=1,\ldots, E_h,\\
  \dfrac{\mathrm du_i}{\mathrm dt}&=
   \frac{g_i^L + {f}_i^*}{m_i}=a_i(\bar u_i^*-u_i),\qquad i=1,\ldots,N_h,
   \label{eq:IDPcritCG}
\end{align}
where $\bar U_e^*\in[U_e^{\min},U_e^{\max}]\subseteq\mathcal G$ and 
$\bar u_i^*\in[u_i^{\min},u_i^{\max}]\subseteq\mathcal G$ are some
locally bound-preserving intermediate states. Then the
semi-discrete scheme  \eqref{eq:EG-MCL} is  LED and
its IDP property can be shown using Theorem 1 from
\cite{kuzmin2022timelim}.
For fully discrete schemes,
preservation of local bounds by forward Euler stages
is guaranteed if the time step $\Delta t$ satisfies the CFL conditions
 \cite{kuzmin2023,kuzmin2022timelim}
\begin{align}
  \Delta tA_e&\le 1,\qquad e=1,\ldots,E_h,\label{eq:CFL1}\\
\Delta ta_i&\le 1,\qquad  i=1,\ldots,N_h.\label{eq:CFL2}
\end{align}

Let us now formulate additional constraints that imply the
validity of semi-discrete entropy inequalities for a particular
entropy pair $\{\eta(u),\mathbf q(u)\}$. According to Tadmor's
theory \cite{tadmor2003entropy}, a cell entropy inequality for $U_e$ holds under the sufficient condition (cf. \cite{KuzminDGentropy,KuHaRu2021,kuzmin2023})
\begin{equation}\label{eq:tadmorDG}
  \left(v(U_{e}) - v(U_{e'})\right)F_{ee'}^*  \leq |S_{ee'}|\left(
  (\boldsymbol{\psi}(U_{e'})-\boldsymbol{\psi}(U_{e}))\cdot\mathbf n_{ee'}
+ \left(v(U_{e}) - v(U_{e'})\right)H_{ee'}^{\mathbb{Q}_0}
  \right),
\end{equation}
where $H_{ee'}^{\mathbb{Q}_0}=\mathcal F(U_e^-,U_e^+;\mathbf n_{ee'})$ is
the low-order LLF flux, $v(u)=\eta'(u)$ is the entropy variable and
$\boldsymbol{\psi}(u)=v(u)\mathbf f(u)-\mathbf q(u)$ is the
entropy potential. In essence, Tadmor's condition \eqref{eq:tadmorDG}
imposes an upper bound
on the rate  $(v(U_{e}) - v(U_{e'}))F_{ee'}^*$ of entropy production
by the flux $F_{ee'}^*$.

The 
rate of entropy production by the element contribution
$f_i^{e,*}$ is given by $(v(u_i)-\bar v_e)f_i^{e,*}$, where
$\bar v_e=\frac{1}{|\mathcal N^e|}\sum_{j\in\mathcal N^e}v(u_j)$ is
the arithmetic mean of $|\mathcal N^e|$ nodal entropy variables.
Adapting Tadmor's theory to the CG setting as in \cite{KuHaRu2021,kuzmin2020g},
we impose the constraint
\begin{equation}\label{eq:tadmorCG}
  (v(u_i)-\bar v_e)f_i^{e,*}\le\sum_{j\in\mathcal N^e\backslash\{i\}}
  |\mathbf c_{ij}^e|
  \left((\boldsymbol{\psi}(u_j)-\boldsymbol{\psi}(u_i))\cdot\mathbf n_{ij}^e
    +   (v(u_i)-v(u_j))H_{ij}^{e,L}\right),
\end{equation}
where
$$
H_{ij}^{e,L}=\frac{\mathbf f(u_j)+\mathbf f(u_i)}2\cdot
\mathbf n_{ij}^e-\frac{d_{ij}^e}{2|\mathbf c_{ij}^e|}(u_j-u_i),
\qquad \mathbf n_{ij}^e=\frac{\mathbf c_{ij}^e}{|\mathbf c_{ij}^e|}.
$$
Note that $\sum_{i\in\mathcal N_e}
(v(u_i)-\bar v_e)f_i^{e,*}=
\sum_{i\in\mathcal N_e}
v(u_i)f_i^{e,*}$ because $\sum_{i\in\mathcal N_e}f_i^{e,*}=0$
by \eqref{eq:zerosumsMCL}. As shown in \cite{kuzmin2020g},
this zero sum property of the element vector $(f_i^{e,*})_{i\in\mathcal N^e}$
also implies the existence of numerical fluxes
$f_{ij}^{e,*}=-f_{ji}^{e,*}$ such that
$f_i^{e,*}=\sum_{j\in\mathcal N^e\backslash\{i\}}f_{ij}^{e,*}$.
It follows that
$$
\sum_{i\in\mathcal N_e}(v(u_i)-\bar v_e)f_i^{e,*}
=\sum_{i\in\mathcal N_e}v(u_i)
\sum_{j\in\mathcal N^e\backslash\{i\}}
f_{ij}^{e,*}
=\sum_{i\in\mathcal N_e}\sum_{j\in\mathcal N^e\backslash\{i\}}
\frac{v(u_i)-v(u_j)}2f_{ij}^{e,*}.
$$
Using this auxiliary result, the entropy stability 
property of the semi-discrete scheme for CG nodal values
can be established as in \cite[Thm.~3]{kuzmin2020g}.
In practice, there is no need to calculate the fluxes 
$f_{ij}^{e,*}$ because the stability condition \eqref{eq:tadmorCG}
is imposed directly on the sum $f_i^{e,*}$.

It remains to select or devise practical algorithms for enforcing the
above constraints. In the next two subsections, we present
the old and new limiters that we use for this purpose.

\subsection{Monolithic limiting for EG cell averages }\label{sec:limiterDG}

The flux limiting framework developed in
\cite{kuzmin2021,KuzminDGentropy} for DG methods
and in \cite{Kuzmin-EG} for an EG discretization of the linear 
advection equation is directly applicable to our nonlinear
problem for cell averages. By definition \eqref{eq:LLF}
of the LLF flux, the equation for $U_e$ can be written as
\begin{align}\notag
  |K_e|\dfrac{\mathrm dU_e}{\mathrm dt}&=
  -\sum_{e'\in\mathcal Z_e}|S_{ee'}|
  \left(\dfrac{\mathbf f(U_e)+\mathbf f(U_{e'})}{2}\cdot
  \mathbf n_{ee'} - \dfrac{\lambda_{ee'}}{2}(U_{e'} - U_e)
  \right)+\sum_{e'\in\mathcal Z_e}F_{ee'}^*\\
  &=\sum_{e'\in\mathcal Z_e}(|S_{ee'}|\lambda_{ee'}
  (\bar U_{ee'}-U_e)+F_{ee'}^*)=\sum_{e'\in\mathcal Z_e}
  |S_{ee'}|\lambda_{ee'}(\bar U_{ee'}^*-U_e),
  \label{eq:barformDG}
\end{align}
where \cite{kuzmin2020,kuzmin2021,kuzmin2022timelim}
\begin{equation*}
  \Bar{U}_{ee'} = \dfrac{U_{e'} + U_e}{2} -  \dfrac{\mathbf f(U_{e'}) - \mathbf f( U_{e})}{2 \lambda_{ee'}}\cdot\mathbf n_{ee'},\qquad
  \Bar{U}_{ee'}^*=\Bar{U}_{ee'}+\frac{F_{ee'}^*}{|S_{ee'}|\lambda_{ee'}}.
\end{equation*}
The low-order \emph{bar state} $\Bar{U}_{ee'}=\Bar{U}_{e'e}$
is a convex combination
of $U_{e'}$ and $U_e$. This can be shown as in \cite{kuzmin2020} using
the mean value theorem and the definition of $\lambda_{ee'}$ as a local
upper bound for the wave speed $|\mathbf f'(u)\cdot\mathbf n_{ee'}|$.
It follows that $\Bar{U}_{ee'}\in\mathcal G=[u^{\min},u^{\max}]$
whenever $U_e,U_{e'}\in\mathcal G$. For any choice of local
bounds $U_e^{\min}\in [u^{\min},\Bar{U}_{ee'}]$ and
$U_e^{\max}\in [\Bar{U}_{ee'},u^{\max}]$, the inequality constraints
 \begin{equation}\label{eq:local-bounds-eg=bp}
   U^{\min}_e \leq \bar U_{ee'}
   +  \dfrac{F_{ee'}^*}{\lambda_{ee'}} \leq U^{\max}_e,
   \qquad  U^{\min}_{e'} \leq \bar U_{ee'}
   -  \dfrac{F_{ee'}^*}{\lambda_{ee'}} \leq U^{\max}_{e'},
 \end{equation}
 can always be enforced by limiting the magnitude of the flux
 $F_{ee'}^*=-F_{e'e}^*$. The validity of \eqref{eq:local-bounds-eg=bp}
 implies that $\bar U_{ee'}^*\in[U^{\min}_e,U_e^{\max}]$
and $\bar U_{e'e}^*\in[U^{\min}_{e'},U_{e'}^{\max}]$. The
IDP and LED properties follow from the fact that \eqref {eq:barformDG}
 can be written in the form \eqref{eq:IDPcritDG} with \cite{kuzmin2022timelim}
 $$
 A_e=\frac{1}{|K_e|}\sum_{e'\in\mathcal Z_e}|S_{ee'}|\lambda_{ee'},
 \qquad \bar U_e^*=\frac{1}{A_e|K_e|}\sum_{e'\in\mathcal Z_e}|S_{ee'}|
 \lambda_{ee'}\bar U_{ee'}^*,
 $$
 where $\bar U_e^*$ is a convex combination of the 
bound-preserving intermediate states $\bar U_{ee'}^*$.

Following Kuzmin et al. \cite{Kuzmin-EG}, we define the local
bounds for EG cell averages as follows:
  \begin{equation*}
    \begin{split}
        U_e^\text{max} &:= \max\left\{\max_{\partial K_e\cap \Gamma^\text{in}} u_\text{in},\max_{i\in \mathcal N^e} \Big\{\max_{e' \in\mathcal E_i} U_{e'}, u_i \Big\}\right\},\\
        U_e^\text{min} &:=\min\left\{\min_{\partial K_e\cap \Gamma^\text{in}} u_\text{in},\min_{i\in \mathcal N^e} \Big\{\min_{e' \in\mathcal E_i} U_{e'}, u_i \Big\}\right\}
    \end{split}
  \end{equation*}
  and calculate bound-preserving (BP) approximations
$F_{ee'}^{\rm BP}=-F_{e'e}^{\rm BP}$ to the target fluxes $F_{ee'}=-F_{e'e}$
 defined by \eqref{eq:antidiffDG}  using the finite volume / DG version 
  {$$
     F_{ee'}^{\rm BP} = 
            \begin{cases}
              \min\left\{F_{ee'},|S_{ee'}|\lambda_{ee'} \min\left\{U^\text{max}_e - \Bar{U}_{ee'}, {\Bar{U}_{e'e} - U^\text{min}_{e'}} \right\}\right\}
         & \mbox{if}\ F_{ee'}>0,\\
         \max\left\{F_{ee'},       
        |S_{ee'}|  \lambda_{ee'}{\max}\left\{U^\text{min}_e - \Bar{U}_{ee'}, {\Bar{U}_{e'e} - U^\text{max}_{e'} }\right\}\right\}
              & \mbox{if}\ F_{ee'}\le 0
            \end{cases}
$$}
of the monolithic convex limiting formula derived in \cite{kuzmin2020}
in the context of CG discretizations.
  
 To enforce Tadmor's condition \eqref{eq:tadmorDG} using a limiter-based
 entropy fix (cf. \cite{KuzminDGentropy,kuzmin2020g,KuHaRu2021}),
 we constrain the final flux
 $F_{ee'}^*=\alpha_{ee'}F_{ee'}^{\rm BP}$ of the MCL scheme
 \eqref{eq:barformDG} using the correction factor 
            \begin{equation*}
    \alpha_{ee'} = \begin{cases}
        \dfrac{Q_{ee'}}{P_{ee'}}\qquad &\text{if }P_{ee'} > Q_{ee'},\\
    1\qquad &\text{otherwise,}
    \end{cases}
\end{equation*}
where 
\begin{equation*}
  P_{ee'} = \left(v(U_{e'}) -v(U_{e}) \right)F_{ee'}^{\rm BP},\qquad
  Q_{ee'}=Q_{ee'}^++\min\left\{0,Q_{ee'}^-\right\}\ge 0.
\end{equation*}
The components $Q_{ee'}^\pm$ of the entropy-dissipative bound $Q_{ee'}$
are defined by (cf. \cite{KuzminDGentropy})
\begin{align*}
  Q_{ee'}^+&=|S_{ee'}|
 \left(v(U_{e'}) - v(U_{e})\right)\dfrac{\lambda_{ee'}}{2} \left(U_{e'} - U_e\right), \\
 Q_{ee'}^- &=|S_{ee'}|\left(
 \boldsymbol{\psi}(U_{e'}) - \boldsymbol{\psi}(U_e) - \left(v(U_{e'}) - v(U_e)\right)
   \frac{\mathbf f(U_{e'}) + \mathbf f(U_{e})}2\right)\cdot \mathbf n_{ee'}.
\end{align*}
The nonnegativity of $Q_{ee'}$ follows from the fact that the low-order
LLF scheme is entropy stable (see, e.g., \cite{kuzmin2023}).
It is easy to verify that $Q_{ee'}$ is bounded
above by the right-hand side of\eqref{eq:tadmorDG}
and that $(v(U_{e'}) - v(U_{e}))F^*_{ee'}\le Q_{ee'}$ for our choice of the
correction factor $\alpha_{ee'}=\alpha_{e'e}$.

Since the above limiting strategy for $F^*_{ee'}$ is not new, we refer the
reader to the original publications \cite{kuzmin2021,Kuzmin-EG,KuzminDGentropy}
and to the book \cite{kuzmin2023}
for detailed explanations and proofs.

\subsection{Monolithic limiting for CG nodal values}\label{sec:limiterCG}

The limiter that we apply to the element contributions $f_i^{e,*}$
is an entropy-stable extension of the element-based MCL algorithms proposed 
in~\cite{Kuzmin-EG},\cite[Sec. 6.3]{kuzmin2023}.
We first express the flux-corrected equation for the CG nodal
value $u_i$ in terms of intermediate states as follows:
\begin{align}\notag
m_i\frac{\mathrm d u_i}{\mathrm d t}
=& -\sum_{e\in\mathcal E_i}\sum_{e'\in\mathcal E_e'}\sigma_{i,ee'}
[\mathcal F(u_i,\hat u_{i,ee'};\mathbf n_{ee'})-\mathbf f(u_i)\cdot\mathbf n_{ee'}]\\ \notag
&+\sum_{e\in\mathcal E_i}\sum_{j\in\mathcal N^e}[d_{ij}^eu_j-\mathbf c_{ij}^e
  \cdot\mathbf{f}(u_j)]+\sum_{e\in\mathcal E_i}f_i^{e,*}
\\ \notag
   =& \sum_{e\in\mathcal E_i}\sum_{e'\in\mathcal E_e'}\frac{\sigma_{i,ee'}}2
   \left[\lambda_{ee'}(\hat u_{i,ee'}-u_i)-
     (\mathbf f(\hat u_{i,ee'})-\mathbf f(u_i))\cdot\mathbf
      n_{ee'}\right]\\
      & +\sum_{e\in\mathcal E_i}\sum_{j\in\mathcal N^e\backslash\{i\}}
      \left[d_{ij}^e(u_j-u_i)-
        (\mathbf{f}(u_j)-\mathbf{f}(u_i))\cdot\mathbf c_{ij}^e\right]
+\sum_{e\in\mathcal E_i}f_i^{e,*}
      \notag \\
      =& \sum_{e\in\mathcal E_i}\left[\sum_{e'\in\mathcal E_e'}
      \sigma_{i,ee'}\lambda_{ee'}(\bar u_{i,ee'}-u_i)
      +\sum_{j\in\mathcal N^e\backslash\{i\}}
      2d_{ij}^e(\bar u_{ij}^e-u_i)+f_i^{e,*}\right]\notag\\
      =&\sum_{e\in\mathcal E_i}[\gamma_i^e(\bar u_i^{e}-u_i)+f_i^{e,*}]
        =\sum_{e\in\mathcal E_i}\gamma_i^e(\bar u_i^{e,*}-u_i),
       \label{eqn:MCLschemeCG}       
\end{align}
where 
\begin{gather*}
  \bar u_{i,ee'} = \dfrac{\hat u_{i,ee'} + u_i}{2}
  - \dfrac{\left(\mathbf f(\hat u_{i,ee'})
    - \mathbf f(u_i) \right)\cdot \mathbf n_{ee'}}{2\lambda_{ee'}},\qquad
  \bar u^e_{ij} = \dfrac{u_j + u_i}{2} - \dfrac{\left(\mathbf f(u_j) - \mathbf f(u_i) \right)\cdot \mathbf c^e_{ij}}{2d^e_{ij}},\\
  \bar u_i^{e}=\frac{1}{\gamma_i^e}\left[\sum_{e'\in\mathcal E_e'}
    \sigma_{i,ee'}\lambda_{ee'}\bar u_{i,ee'}
    +\sum_{j\in\mathcal N^e\backslash\{i\}}
    2d_{ij}^e\bar u^e_{ij}\right],\qquad
  \bar u_i^{e,*}=\bar u_i^{e}+\frac{f_i^{e,*}}{\gamma_i^e},
  \\
    \gamma_i^e=\sum_{e'\in\mathcal E_e'}
    \sigma_{i,ee'}\lambda_{ee'}
    +\sum_{j\in\mathcal N^e\backslash\{i\}}
    2d_{ij}^e.
  \end{gather*}
 Therefore, equation \eqref{eqn:MCLschemeCG} can be written in the form
 \eqref{eq:IDPcritCG} with
 $$
 a_i=\frac{1}{m_i}\sum_{e\in\mathcal E_i}\gamma_i^e,\qquad
 \bar u_i^*=\frac{1}{a_im_i}\sum_{e\in\mathcal E_i}
 \gamma_i^e\bar u_i^{e,*}.
 $$
To ensure that the intermediate state $\bar u_i^{e,*}$ is bounded by


$$
u_i^{\max} = \max\left\{
                u_i, \max_{e\in\mathcal E_i}\bar{u}_i^e \right\} \  \text{ and } \
u_i^{\min} = \min\left\{
               u_i, \min_{e\in\mathcal E_i}\bar{u}_i^e \right\},
$$
 we impose the MCL constraints
 \begin{equation}\label{eq:local-bounds-BP}
   u_i^{\min} \leq \bar u^e_i +
   \dfrac{f^{e,*}_i}{\gamma^e_{i}}
   \leq u^{\max}_i.
 \end{equation}
 To satisfy the entropy stability condition \eqref{eq:tadmorCG},
 we additionally require that
\begin{equation}\label{eq:local-bounds-ES}
  (v(u_i) - \bar v_e ) f^{e,*}_i \leq \sum_{j\in\mathcal N^e
    \backslash\{i\}}[Q_{ij}^{e,+}+
  \min\{0,Q_{ij}^{e,-}\}]=:Q^{e}_i,
\end{equation}
where $Q_i^e$ is a nonnegative production bound depending on
\begin{align*}
  Q^{e,+}_{ij} &=
  (v(u_j) - v(u_i)) \dfrac{d^e_{ij}}{2}(u_j - u_i),\\
      Q_{ij}^{e,-} &= 
        \left(\bm \psi(u_j) - \bm \psi(u_i) - (v(u_j) - v(u_i))
        \dfrac{\mathbf f(u_j) + \mathbf f(u_i)}{2}\right)\cdot \mathbf c^e_{ij}.
\end{align*}
        For $f_i^{e,*}=\alpha_i^ef_i^e$
        with $\alpha_i^e\in[0,1]$,
condition \eqref{eq:local-bounds-ES} is equivalent to
$$
\alpha_i^eP_i^e\le Q_i^e,\qquad P_i^e:= (v(u_i) - \bar v_e ) f^{e}_i. 
$$

Adopting a clip-and-scale limiting strategy (cf.
\cite{lohmann2017,Kuzmin-EG,kuzmin2023}), we first apply
$$
\qquad    \alpha_i^e= \begin{cases}
            \dfrac{Q^{e}_i}{P^e_i}&\qquad \text{if }P^e_i > Q^{e}_i,\\
            1&\qquad \text{otherwise}
\end{cases}
$$
to  $f_i^e$ defined by \eqref{eq:antidiffCG}.
Then we calculate the clipped element contributions
$$
\tilde f_i^{e,*}=\begin{cases}
\min\{\alpha_i^ef_i^e,\gamma_i^e(u_i^{\max}-\bar u_i^e)\} &\mbox{if}\ f_i^e>0,\\
\max\{\alpha_i^ef_i^e,\gamma_i^e(u_i^{\min}-\bar u_i^e)\} &\mbox{if}\ f_i^e\le 0.
\end{cases}
    $$

        It is easy to verify that conditions \eqref{eq:local-bounds-BP} and
        \eqref{eq:local-bounds-ES} hold for 
        $f_i^{e,*}=\tilde f_i^{e,*}$. However, this definition of $f_i^{e,*}$
      does not generally ensure the validity of the zero sum condition
        $
\tilde f_e^++\tilde f_e^-=0
$
for
$$\tilde f_e^+=\sum_{i\in\mathcal N^e}\max\{0,f_i^{e,*}\},\qquad
\tilde f_e^-=\sum_{i\in\mathcal N^e}\min\{0,f_i^{e,*}\}.
$$
In the scaling stage of our entropic MCL procedure, we calculate (cf.~\cite{lohmann2017,Kuzmin-EG})
\begin{equation}
  f_i^{e,*}=\begin{cases}
  \left(-\frac{\tilde f_e^-}{\tilde f_e^+}\right)\tilde f_i^{e,*}
  & \mbox{if}\  \tilde f_e^++\tilde f_e^->0 \ \mbox{and}\ \tilde f_i^{e,*}>0,\\
   \left(-\frac{\tilde f_e^+}{\tilde f_e^-}\right)\tilde f_i^{e,*}
   & \mbox{if}\  \tilde f_e^++\tilde f_e^-<0 \ \mbox{and}\ \tilde f_i^{e,*}<0,\\
\tilde f_i^{e,*} & \mbox{otherwise}.
  \end{cases}
\end{equation}
This final result satisfies the zero sum condition
$\sum_{i\in\mathcal N^e}f_i^{e,*}=0$ in addition to \eqref{eq:local-bounds-BP} and
        \eqref{eq:local-bounds-ES}. 
        


\begin{remark}
\label{rem:1}
By default, the monolithic limiting strategy for (\ref{eq:EG-MCL}), as presented in Sections \ref{sec:limiterDG} and \ref{sec:limiterCG},  enforces both preservation of local bounds and entropy stability. The inequality constraints associated with either of these properties can be modified or deactivated for testing purposes. For example, flux limiters based solely on conditions (\ref{eq:local-bounds-eg=bp}) and (\ref{eq:local-bounds-BP}) would disregard the entropy stability conditions, as in \cite{kuzmin2023, Kuzmin-EG, kuzmin2021}. Several comparisons between bound-preserving limiting strategies without and with optional entropy fixes are performed in Section~\ref{sec:4:Numerical Results}.

\end{remark}

        \subsection{Constrained projection of output data} \label{sec:FCR}

        The limiting techniques presented so far guarantee preservation of local bounds for the EG cell averages and CG nodal values. However, no discrete maximum principle generally holds for the restriction $u_h^\text{EG}|_{K_e} = u_h|_{K_e} + \delta u_e = u_h|_{K_e} + \left(U_e - \bar u_e\right)$ of the EG solution to a cell $K_e$. When it comes to visualizing the results or using them as input data in solvers for other equations, we project $u_h^\text{EG}$ into the CG space $V_h^{\rm CG}$ using the flux-corrected remapping (FCR) algorithm presented in \cite{Kuzmin-EG}. The nodal values of the standard $L^2$ projection $u_h^H= \sum_{j=1}^{N_h} u_j^H\varphi_j$ satisfy
        $$
        \sum_{j\in\mathcal N_i}m_{ij}u_j^H=\sum_{e\in\mathcal E_i}\int_{K_e}u_h^{\EG}
        \mathrm d\mathbf x,\qquad i=1,\ldots,N_h.
        $$
  The FCR scheme yields a bound-preserving approximation $u_h^\text{FCR} = \sum_{j=1}^{N_h} u_j^\text{FCR}\varphi_j$ such that
\begin{equation}\label{eq:FCR-local}
  u^{\min,\rm FCR}_i \leq u_i^\text{FCR}
=u_i^L+\frac{1}{m_i}\sum_{e\in\mathcal E_i}\alpha_e^{\rm FCR}f_i^{e,\rm FCR}
  \leq u^{\max,\rm FCR}_i ,\qquad  i =1,\ldots, N_h,
\end{equation}
where 
$$
u_i^{\max,\rm FCR}=\max\left\{\max_{e\in\mathcal E_i} U_e, \max_{j\in\mathcal N_i} u_j\right\},\qquad
u_i^{\min,\rm FCR}=
 \min\left\{\min_{e\in\mathcal E_i} U_e, \min_{j\in\mathcal N_i} u_j\right\},
 $$
 $$
 u_i^L= \frac{1}{m_i}\sum_{e\in\mathcal E_i} m^e_i U_e,\qquad
 f_i^{e,\rm FCR}=
 \int_{K_e} \varphi_i \left(u^{H}_i -
 u^{H}_h + u^\text{EG}_h - U_e
    \right)\mathrm d\mathbf x,
    $$
\begin{equation*}
        \alpha_e^{\rm FCR} = \min_{i\in\mathcal N^e} \begin{cases}
            \min\left\{ 1, \dfrac{m^e_i (u_i^{\max,\rm FCR} - u^{L}_i)}{f^{e,\rm FCR}_i}\right\}\qquad &\text{if }f^{e,\rm FCR}_i >0,\\
            \min\left\{ 1, \dfrac{m^e_i (u_i^{\min,\rm FCR} - u^{L}_i)}{f^{e,\rm FCR}_i}\right\}\qquad &\text{if }f^{e,\rm FCR}_i <0,\\
            1\qquad &\text{otherwise}.
        \end{cases}
    \end{equation*}
 We plot $u_h^\text{FCR}$ instead of $u_h^\text{EG}$ in the figures that we 
present in Section \ref{sec:4:Numerical Results}. However, we do not overwrite
 $u_h^\text{EG}$ by $u_h^\text{FCR}$ and calculate the error
 norms using $u_h^\text{EG}$ rather than $u_h^\text{FCR}$.

\section{Temporal discretization}\label{sec:SSp-RK}

We discretize system \eqref{eq:EG-MCL} in time using Heun's method, a second-order explicit strong stability preserving Runge--Kutta  method with two forward Euler updates (SSP-RK2). Let $\Delta t=T/N_T$ denote the constant time step corresponding to a uniform subdivision of the time interval $[0,T]$ into $N_T\in\mathbb N$ subintervals. Individual components of the global vectors
$$
\mathbf U^{n} = \left(U^{n}_e\right)_{e = 1}^{E_h}\quad \mbox{and}\quad
    \quad\mathbf u^{n} = \left(u^{n}_i\right)_{i=1}^{N_h}
    $$
    represent approximations to the EG cell averages and CG nodal values,
    respectively, at
    the time level $t^n = n\Delta t,\  0\leq n\leq N_T$. Applying 
    Heun's method to the coupled subproblems \eqref{eq:IDPcritDG} and
    \eqref{eq:IDPcritCG}, we advance the degrees of freedom
    $U^{n}_e$ and $ u^{n}_i$
    to the time level $t^{n+1}$
    as follows:
\begin{itemize}
    \item[1.]First forward Euler step
    \begin{alignat*}{3}
      U^{(n,1)}_e &= U^n_e + \Delta t A_e^n(\bar U_e^{n,*}-U_e^n), &&\qquad e = 1,\ldots, E_h,\\
    u^{(n,1)}_i &= u^n_i + \Delta ta_i^n(\bar u_i^{n,*}-u_i^n), &&\qquad i = 1,\ldots, N_h.
\end{alignat*}
\item[2.]Second forward Euler step
    \begin{alignat*}{3}
      U^{(n,2)}_e &= U^{(n,1)}_e + \Delta t A_e^{(n,1)}(\bar U_e^{(n,1),*}-U_e^{(n,1)}), &&\qquad e = 1,\ldots, E_h,\\
      u^{(n,2)}_i &= u^{(n,1)}_i + \Delta t a_i^{(n,1)}(\bar u_i^{(n,1),*}-u_i^{(n,1)}), &&\qquad i = 1,\ldots, N_h.
    \end{alignat*}
\item[3.] Final SSP-RK stage:
\begin{alignat*}{3}
    U^{n+1}_e &= \frac{1}{2}\left(U^{(2,n)}_e + U^{n}_e\right), &&\qquad e = 1,\ldots, E_h,\\
    u^{n+1}_i &= \frac{1}{2}\left(u^{(2,n)}_i + u^{n}_i\right), &&\qquad i = 1,\ldots, N_h.
\end{alignat*}
\end{itemize}
If the CFL conditions \eqref{eq:CFL1} and \eqref{eq:CFL2} hold, each step
produces a convex combination of old values. Therefore, discrete maximum
principles are satisfied for the new ones. The DG components
$$\delta u_e^{n+1}=U_e^{n+1}-\bar u_e^{n+1},\qquad e=1,\ldots,E_h$$
are treated as derived quantities in our implementation of
the EG method. The optional FCR postprocessing 
(see Section \ref{sec:FCR}) ensures preservation of
local bounds for projected output data.
      
\section{Numerical results}
\label{sec:4:Numerical Results}

  In this section, we illustrate the capabilities of (individual components of) the proposed algorithm by running several numerical experiments. In particular, we perform grid convergence tests and apply the methods under investigation to two-dimensional benchmark problems. Throughout this section, the acronyms \textbf{HO} (High Order) and \textbf{LO} (Low Order) refer to the EG schemes defined in Sections~\ref{sec.numerics} and~\ref{sec:algrbraic}, respectively. The one denoted by \textbf{BP-ES} enforces both preservation of local bounds and entropy stability, as described in Section \ref{sec:MCL}. The method labeled \textbf{BP} differs from \textbf{BP-ES} in that the entropy fixes are deactivated to check if they are needed to ensure convergence towards correct weak solutions. In all numerical experiments, we use
uniform rectangular meshes and codes implemented using the deal.II finite element library \cite{dealII94}.

\subsection{Example 1. Convergence test for a linear advection equation}\label{example1}
{
First, we run the grid convergence test for the following linear advection equation
\begin{equation}
    \dfrac{\partial u}{\partial t} + \nabla \cdot (\mathbf v  u) = 0,
\end{equation}
where $\mathbf v\equiv (1,0)$ is the constant velocity. The computational domain is given by
$\Omega = (0,1)^2$.
The inflow boundary condition is imposed on $\Gamma^\text{in} = \{(x,y) \in\partial \Omega : x = 0 \}$, while the outflow boundary condition is applied otherwise.
The exact solution is given by 
\begin{equation*}
    u(x,y,t) = \cos(\pi(x-t)),
\end{equation*}
with the initial condition $u_0(x,y)  = \cos(\pi x)$. We choose the square entropy $\eta(u) = u^2/2$ corresponding to the {entropy flux $\mathbf q(u) = \left({u^2}/{2},0\right)$}. This setup extends
 a one-dimensional advection problem with the exact solution 
$u(x,t) = \cos(\pi(x-t))$ into the two-dimensional domain.

We test our LO, HO, BP and BP-ES schemes by running numerical simulations on seven uniformly refined meshes up to the final time ${T = 0.5}$. The mesh sizes range from $h = 2^{-2}$ to $h = 2^{-8}$. For numerical integration in time, we use the SSP-RK2 method {(as presented in Section~\ref{sec:SSp-RK})}. To satisfy the CFL-like conditions \eqref{eq:CFL1} and \eqref{eq:CFL2}, the time step $\Delta t = 0.025$ is also refined, maintaining a constant ratio of $\Delta t / h = 0.1$.
We then use the global norms
\begin{subequations}
\begin{alignat*}{2}
    \|u_h^\text{EG} - u\|_{l^\infty(L^1)} &:= \max_{0\le n=0\le N} \|u_h^{\text{EG},n} - u^n \|_{L^1} \qquad \text{and }\quad 
    \|u_h^\text{EG} - u\|_{l^2(H^1)} &:= \sqrt{\sum_{n=1}^{N} \Delta t\|u_h^{\text{EG},n} - u^n \|^2_{H^1}}
\end{alignat*}
\end{subequations}
to measure the numerical error $u_h^\text{EG}-u$.
The convergence results are presented in Figure \ref{fig:Example 1}.

Figure \ref{fig:Example1-L2-error} shows that our LO scheme achieves the expected first-order convergence rate in the $l^\infty(L^1)$-norm. The BP-ES, BP and HO results exhibit nearly second-order convergence. In the final cycle, the convergence rate is 1.91 for both BP-ES and BP, while the HO error shrinks at the rate 1.98. These results align  with those reported in \cite{KUZMIN2022114428,Jean-Luc-first-order,Jean-Luc-second-order}. The $l^2(H^1)$ errors are shown in Figure \ref{fig:Example-H1-error}. The LO convergence rate is 0.59 in the final cycle, while BP-ES, BP and HO achieve optimal first-order convergence rates 0.98, 0.99 and 1.02, respectively.

\begin{figure}[H]
    \centering
    \begin{subfigure}[t]{0.4\textwidth}
        \centering
        \hspace*{-.7cm}\includegraphics[width=1.2\textwidth]{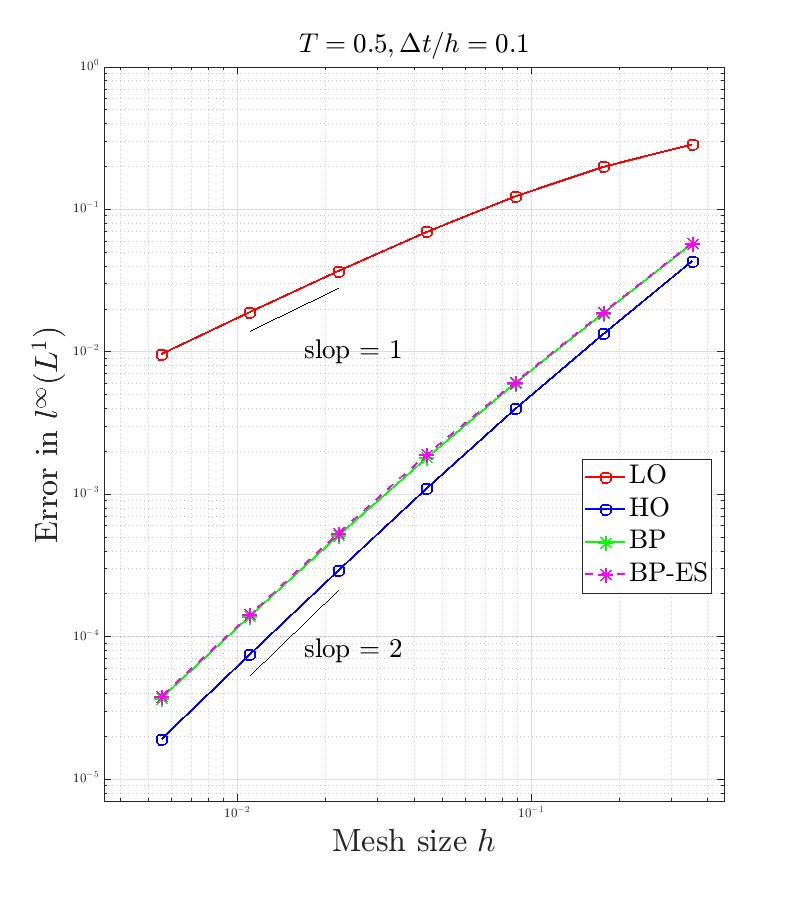}
    
    \vspace*{-0.4cm}\caption{$l^\infty(L^1)$}
    \label{fig:Example1-L2-error}
    \end{subfigure}
    \hfill
    \begin{subfigure}[t]{0.4\textwidth}
        \centering
        \hspace*{-0.7cm}\includegraphics[width=1.2\textwidth]{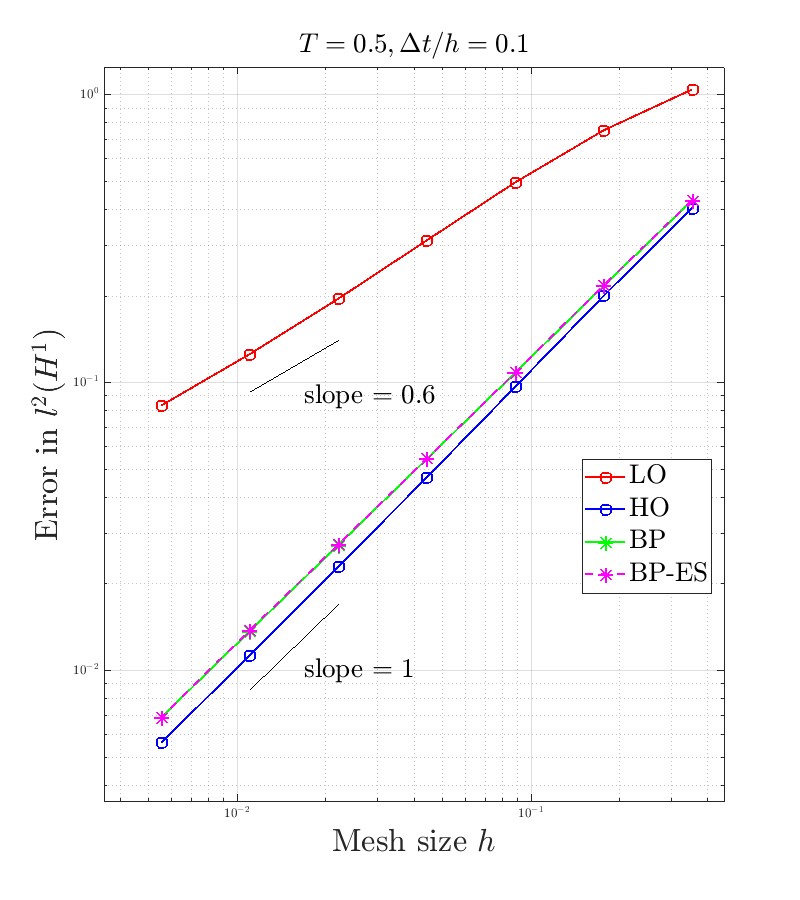}
        
        \vspace*{-0.4cm}\caption{$l^2(H^1)$}
    \label{fig:Example-H1-error}
    \end{subfigure}
    \caption{Example 1. The convergence behavior of EG errors for LO, HO, and BP schemes.}
    \label{fig:Example 1}
\end{figure}

As expected, we obtain the optimal convergence rate for both BP and BP-ES with the employed limiters.
Since the exact solution is smooth, the differences between the constrained EG approximations and the underlying HO scheme are insignificant.

\subsection{Example 2. Convergence test for an inviscid Burgers equation}
\label{sec:Example2}
{
  In this example, we
  examine the performance of our numerical schemes applied to the nonlinear inviscid Burgers equation
\begin{equation}\label{eq:1d-burger}
    \dfrac{\partial u}{\partial t}  +  \nabla \cdot \left(\dfrac{u^2}{2}\mathbf v\right)= 0 
\end{equation}
with $\mathbf v\equiv (1,0)$ in the computational domain $\Omega = (0,1)^2$. The inflow boundary condition is imposed on $\Gamma^\text{in} = \{(x,y)\in\partial \Omega : x  =0\}$, and the outflow boundary condition is applied elsewhere. The initial condition is given by
$u_0(x,y) = \sin(2\pi x)$. For this nonlinear test, we choose the entropy $\eta(u) = u^4/4$ with the corresponding  entropy flux  
$\mathbf q(u) = \left({u^5}/{5},0\right)$.

Note that the unique entropy solution develops a shock at the critical time $t_c = \frac{1}{2\pi}$. However, before $t_c$, the solution remains smooth and can be obtained using the method of characteristics \cite{evans10}. To determine the exact solution value $u(\mathbf x,t)  = u_0(\mathbf x_0)$ for $t<t_c$, we solve the
nonlinear equation $\mathbf x_0 = \mathbf x - u_0(\mathbf x_0)\mathbf v t$
for $\mathbf x_0$ using fixed-point iterations. Similarly to Example \ref{example1}, the problem that we consider in this nonlinear test is essentially one dimensional, although computations are performed in a two-dimensional domain. 

We first test the convergence behavior of the LO, HO, BP, and BP-ES schemes by solving the Burgers equation \eqref{eq:1d-burger} up to the final time $T = 0.1 { < t_c}$ on seven uniformly refined meshes ($h = 2^{-1}$ to $h = 2^{-7}$). The time step  $\Delta t=0.002 h$ is also refined in each cycle.

\begin{figure}[!h]
    \centering
    \begin{subfigure}[t]{0.4\textwidth}
        \centering
        \includegraphics[width=1.1\textwidth,trim=0 0 0 60, clip]{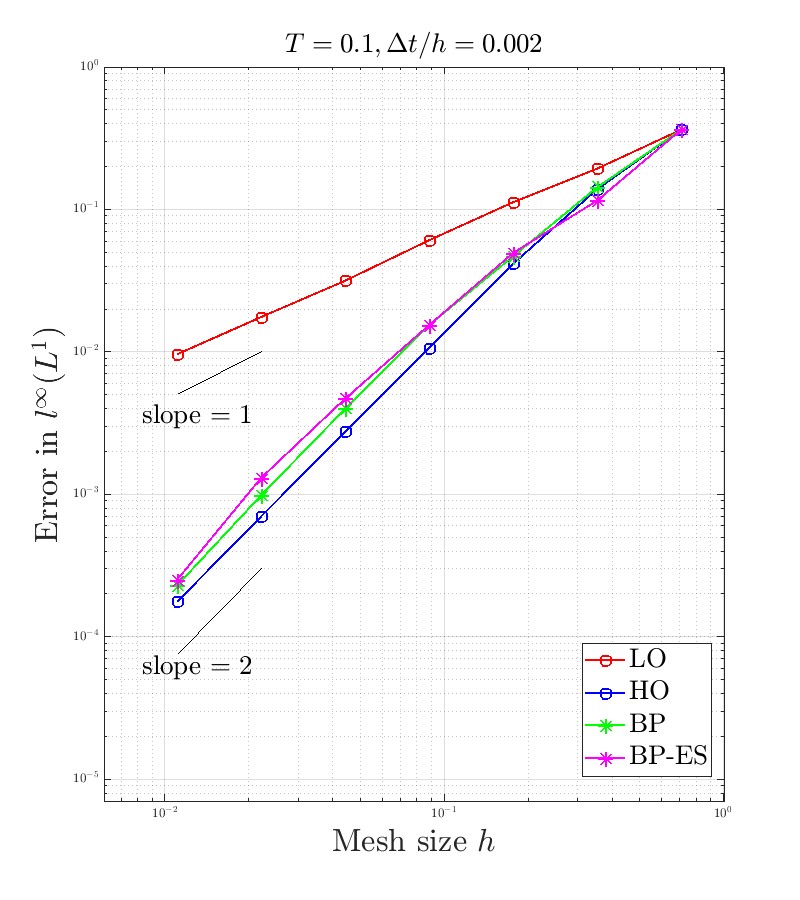}
        \caption{$T = 0.1,\ \Delta t/h=0.002$}
        \label{fig:Example2-Convergence}
    \end{subfigure}\hfill
    \begin{subfigure}[t]{0.46\textwidth}
      \centering
      
        \includegraphics[width=1.1\textwidth]{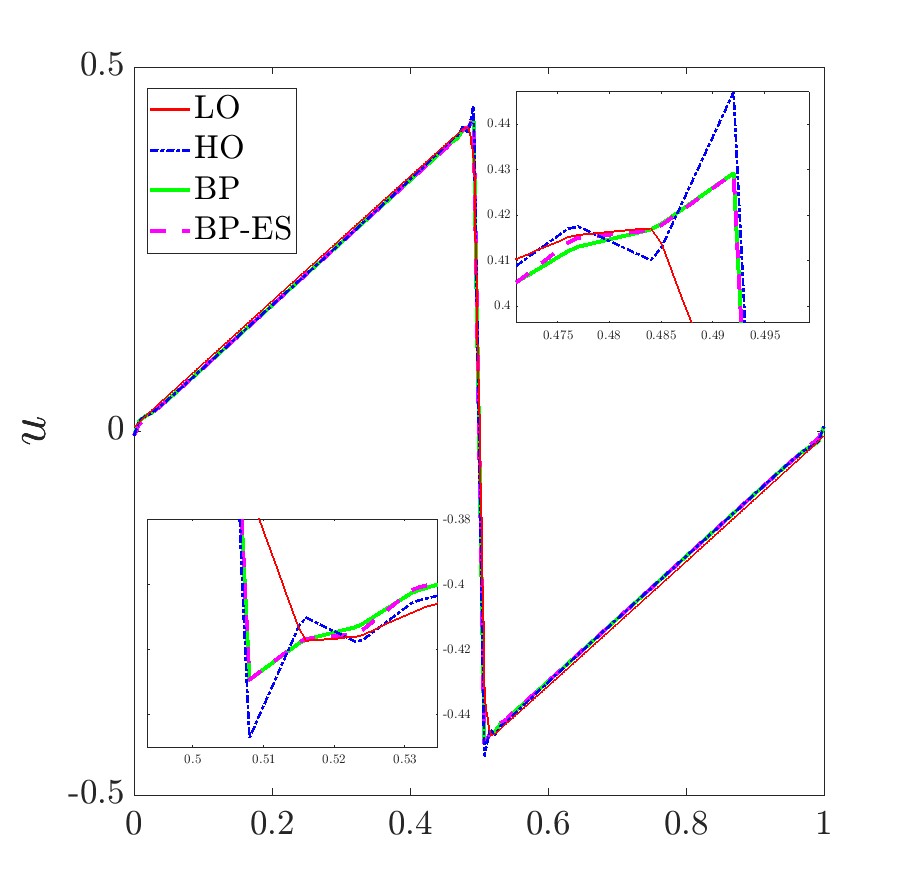}
        \caption{$T = 1.0,\ h=2^{-7},\ \Delta t=0.001$}
        \label{fig:Example2-Shock}
    \end{subfigure}
    \caption{Example 2. The diagrams show (a) error behavior and (b)
      solution profiles along the line $y=0.5$. }
\end{figure}

The $l^\infty(L^1)$-error of the EG solution $u_h^\text{EG}$ is reported in Figure \ref{fig:Example2-Convergence}. The unconstrained HO scheme exhibits second-order convergence at the rate 1.99 in the final cycle. This demonstrates that our baseline EG discretization of the nonlinear problem is second-order accurate, despite the use of the linear Taylor approximation (\ref{eqn:taylor}) and of the group finite element formulation (\ref{eq:GFE}). The  convergence rate of the LO scheme is 0.87 in the final cycle, while the BP and BP-ES convergence rates become as high as 2.10 and 2.38, respectively. This indicates that the first-order error of the LO approximation is fully compensated in the limiting stage.

Next, we extend the final time to $T = 1.0>t_c$ and run a BP-ES simulation of the propagating shock. In this test, we use the mesh with spacing $h = 2^{-7}$ and the time step $\Delta t = 0.001$. For comparison purposes, we also present the LO, HO, and BP results. The profiles of numerical solutions along the middle line $\{y = 0.5\}$ are plotted in Figure \ref{fig:Example2-Shock}.

As expected, the results obtained with the entropy stable and bound-preserving LO scheme are free from spurious oscillations and capture the shock location correctly. The HO scheme is neither entropy stable nor bound preserving. In this example, it produces some oscillations around the shock. It is worth noting that these oscillations are less significant than those observed in simulations with classical CG-$\mathbb Q_1$ methods (cf. Figure 2 in \cite{KUZMIN2023112153}). The stabilizing effect of the EG enrichment $\delta u_h$ was also observed in \cite{Kuzmin-EG} in the context of linear advection problems.

\subsection{Example 3. Convergence test for the two-dimensional KPP problem}

In the third example, we consider a smooth version of the two-dimensional Kurganov-Petrova-Popov (KPP) test \cite{doi:10.1137/040614189}. The nonlinear hyperbolic conservation law
\begin{equation}\label{eq:Example3-KPP}
    \dfrac{\partial u}{\partial t} + \nabla \cdot \mathbf f(u) = 0, \qquad \mathbf f(u) =  \big(\sin(u), \cos(u)\big)
\end{equation}
is solved in the computational domain $\Omega = (-2,2)\times(-2.5, 1.5)$. In our grid convergence study, we use the smooth initial condition
\begin{equation}\label{eq:KPP-smooth-initial}
    u_0(x,y) = \begin{cases}
        \dfrac{\pi}{4}\Big(1 + \dfrac{1}{20}\big(1+\cos(\pi\sqrt{x^2+y^2}\big)
        \Big)\qquad &\text{if }\sqrt{x^2+y^2}\leq 1,\\\\
        \dfrac{\pi}{4}\qquad &\text{otherwise.}
    \end{cases}
\end{equation}
We select the square entropy $\eta (u) = u^2 / 2$, which is 
is associated with the entropy flux \cite{KUZMIN2022114428}
\begin{equation}\label{eq:Example3qu}
 \mathbf q(u) = \left(\cos(u)+ u \sin(u), -\sin(u) + u \cos(u)\right).   
\end{equation}
When it comes to calculating the artificial viscosity coefficients $d^e_{ij}$ and the LLF fluxes using formulas (\ref{eq:d_e_ij}) and (\ref{eq:LLF}), respectively, we use $\lambda = 1$ as a global upper bound for the maximum wave speed. More accurate estimates of that speed can be found in \cite{Jean-Luc-second-order}.

In this test, no analytical solution is available. Therefore, the experimental order of convergence (EOC) for $h = 2^{-7}$ is determined by the differences in the $L^1$- and $L^2-$norms of numerical solutions on three successively refined grids. Given the numerical solutions calculated on meshes with spacings $h$, $2h$, and $4h$, the three-level EOCs are calculated as follows:
\begin{equation*}
    \text{EOC}_{L^1} = \log_2 \left( \dfrac{\|u_{4h} - u_{2h}\|_{L^1}}{\|u_{2h} - u_h\|_{L^1}}\right),\qquad 
    \text{EOC}_{L^2} = \log_2 \left( \dfrac{\|u_{4h} - u_{2h}\|_{L^2}}{\|u_{2h} - u_h\|_{L^2}}\right).
\end{equation*}

 The time step is set to keep the ratio $\Delta t/h = 0.256$ fixed. The values of $\text{EOC}_{L^1}$ and $\text{EOC}_{L^2}$ are output at the final time $T = 1.0$. Table \ref{tab:my_label} presents the results of our grid convergence study. It can be seen that the convergence rates of the HO, BP and BP-ES schemes are approximately twice as high as the EOC of the LO approximation.
 

\begin{table}[H]
    \centering
    \begin{tabular}{c |cccc}
           & LO  &  HO & BP & BP-ES\\
           \hline
     EOC$_{L^1}$ & 0.87 & 1.93 & 1.85 & 1.82\\
     EOC$_{L^2}$ & 0.83 & 1.73 & 1.72 & 1.69
    \end{tabular}
    \caption{Example 3: KPP convergence test with smooth initial conditions, $L^1$ and $L^2$ convergence rates. }
    \label{tab:my_label}
\end{table}

\subsection{Example 4. Two-dimensional rotational KPP test}

{In the final example, we solve the nonlinear conservation law \ref{eq:Example3-KPP} in the same computational domain, $\Omega = (-2,2)\times(-2.5, 1.5)$, but using the discontinuous initial condition \cite{doi:10.1137/040614189}
\begin{equation}\label{eq:KPP-initial-discontinuous}
    u_0(x,y) = \begin{cases}
       \dfrac{7\pi}{2}\qquad &\text{if }\sqrt{x^2 + y^2}\leq 1,\\\\
       \dfrac{\pi}{4}\qquad &\text{otherwise.}
    \end{cases}
\end{equation}
The so-defined classical KPP problem admits infinitely many weak solutions. The unique entropy solution at the final time $T=1.0$ exhibits a two-dimensional rotational wave structure~\cite{Jean-Luc-second-order,doi:10.1137/040614189,KUZMIN2022114428}. The challenge of this test is to ensure that a numerical scheme provides entropy stability to prevent convergence to incorrect weak solutions. In our BP-ES method, we perform algebraic entropy fixes using the square entropy $\eta(u) = u^2/2$ and the entropy flux $\mathbf q(u)$ given by (\ref{eq:Example3qu}).

We run computations on a uniform mesh using the discretization parameters $h= 2^{-7}$ and $\Delta t = 0.001$.
All the other settings are the same as in the Example 3. 

In Figures \ref{fig:Example4-2d-discrete-profiles} and \ref{fig:Example4-3d-profiles}, we compare the EG solutions produced by the LO, BP, and BP-ES schemes at the final time. In addition, we present the BP-ES result obtained using a flux-corrected transport (FCT) algorithm from \cite{Kuzmin-EG} instead of MCL to enforce preservation of local bounds. Details of the FCT limiting procedure can be found in \cite{kuzmin2023,Kuzmin-EG,kuzmin2021,KuzminDGentropy}.

\begin{figure}[h!]
    \centering
    \begin{subfigure}[t]{0.45\textwidth}
        \centering
        \includegraphics[trim={2cm 5cm 0 3.5cm},clip,width =\textwidth]{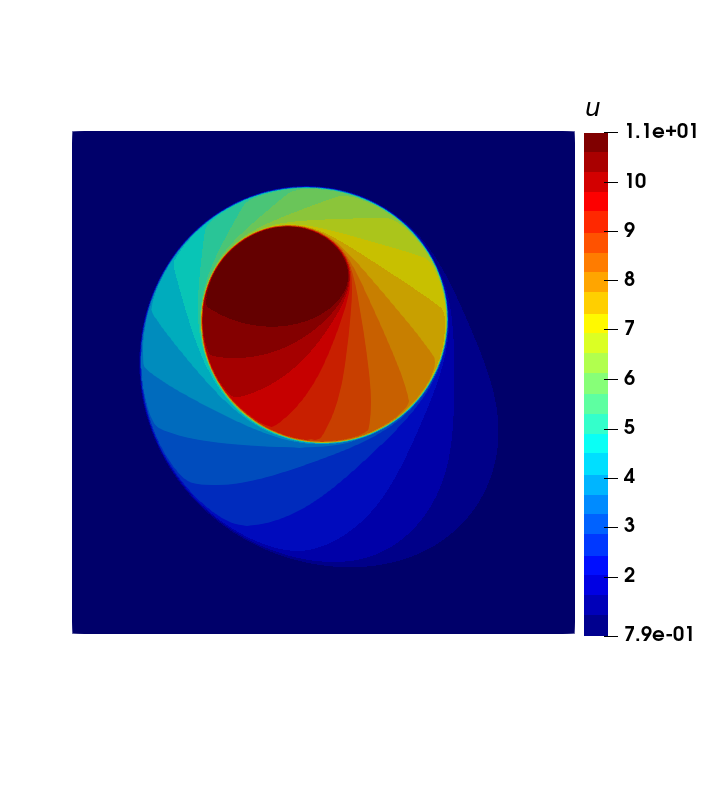}
        \caption{LO}
        \label{fig:Example4-LO}
    \end{subfigure}\hfill
    \begin{subfigure}[t]{0.45\textwidth}
        \centering
        \includegraphics[trim={2cm 5cm 0 3.5cm},clip,width=\textwidth]{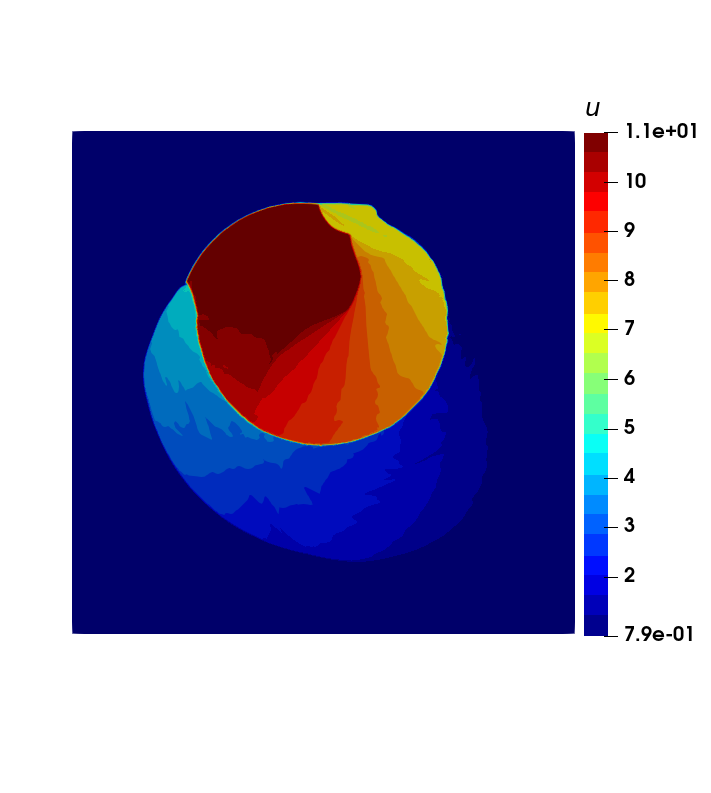}
        \caption{BP}
        \label{fig:Example4-BP}
    \end{subfigure}
    
    \begin{subfigure}[t]{0.45\textwidth}
        \centering
        \includegraphics[trim={2cm 5cm 0 3.5cm},clip,width=\textwidth]{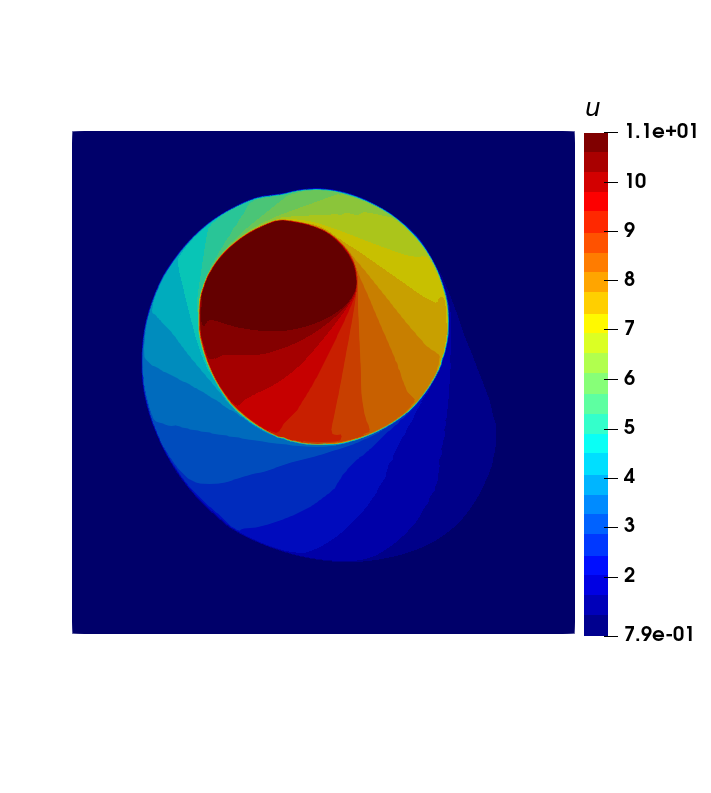}
        \caption{FCT-BP-ES}
        \label{fig:Example4-FCTES}
    \end{subfigure}\hfill
    \begin{subfigure}[t]{0.45\textwidth}
        \centering
        \includegraphics[trim={2cm 5cm 0 3.5cm},clip,width=\textwidth]{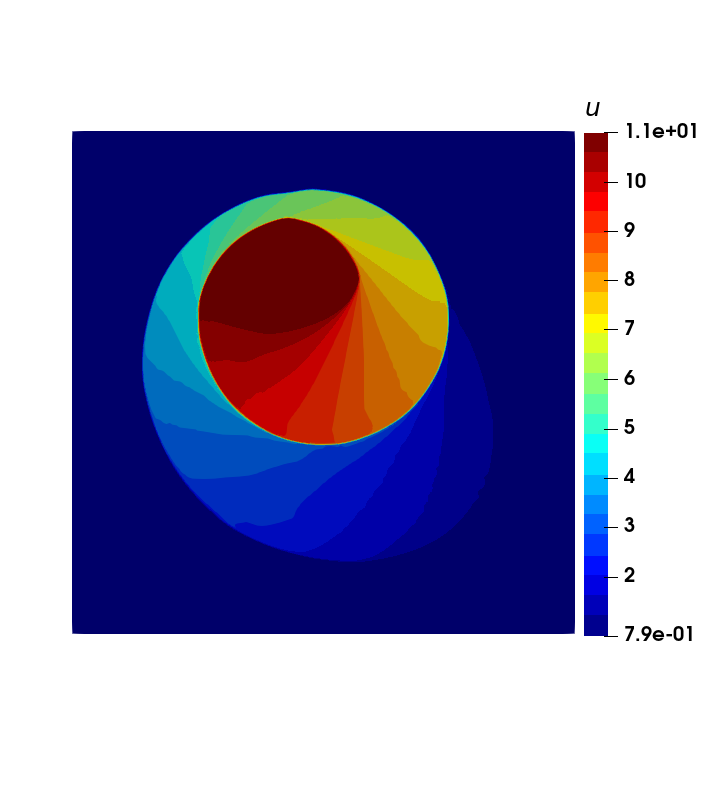}
        \caption{BP-ES}
        \label{fig:Example4-MCL-ES}
    \end{subfigure}
    \caption{Example 4. Numerical solutions produced by different schemes at $T = 1.0$. The colormap of the filled contour plots  represents 25 discrete values in the range $[\frac{\pi}{4},\frac{7\pi}{2}]$. }
    \label{fig:Example4-2d-discrete-profiles}
\end{figure}

\begin{figure}[h!]
    \centering
    \begin{subfigure}[t]{0.4\textwidth}
        \centering
        \includegraphics[width =1.3\textwidth]{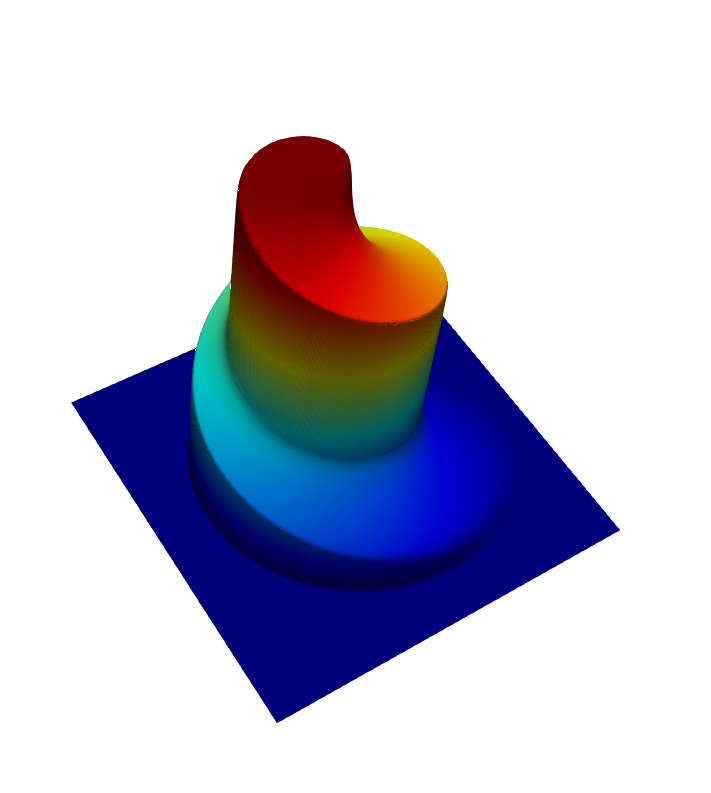}
        \caption{LO}
        \label{fig:Example4-LO-3d}
    \end{subfigure}\hfill
    \begin{subfigure}[t]{0.4\textwidth}
        \centering
        \includegraphics[width=1.3\textwidth]{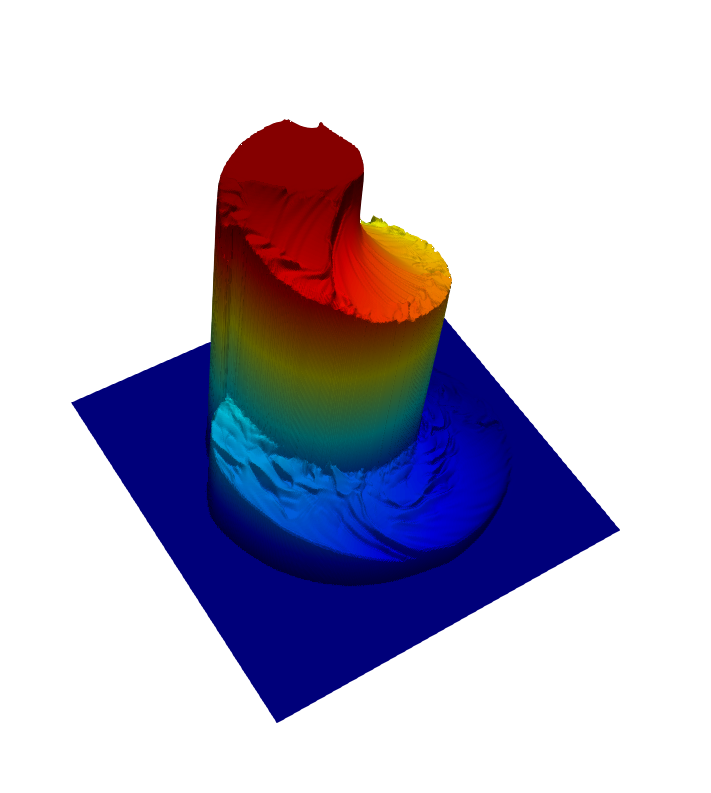}
        \caption{BP }
        \label{fig:Example4-BP-3d}
    \end{subfigure}

    \begin{subfigure}[t]{0.4\textwidth}
        \centering\includegraphics[width=1.3\textwidth]{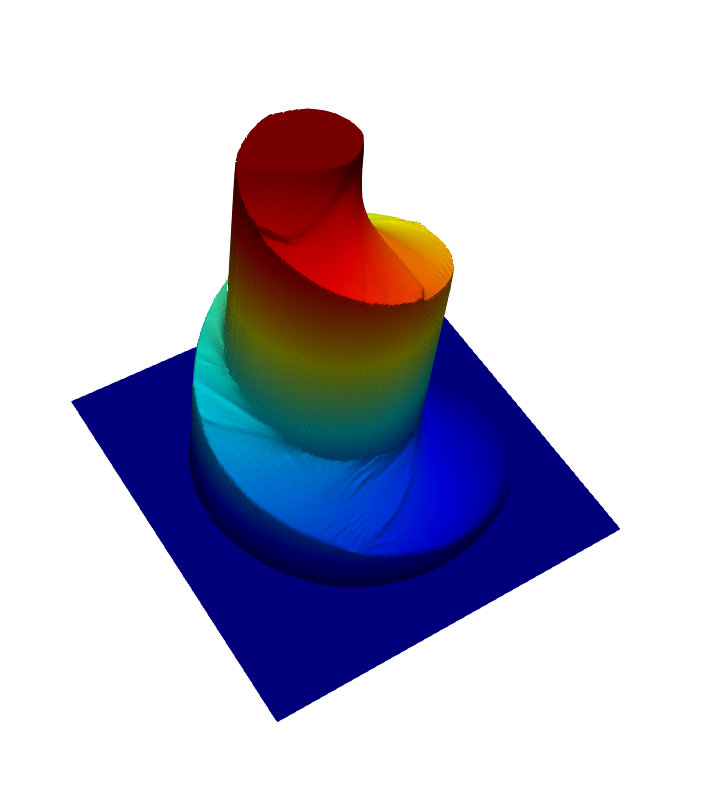}
        \caption{FCT-BP-ES}
        \label{fig:Example4-FCTES-3d}
    \end{subfigure}\hfill
    \begin{subfigure}[t]{0.4\textwidth}
        \centering
        \includegraphics[width=1.3\textwidth]{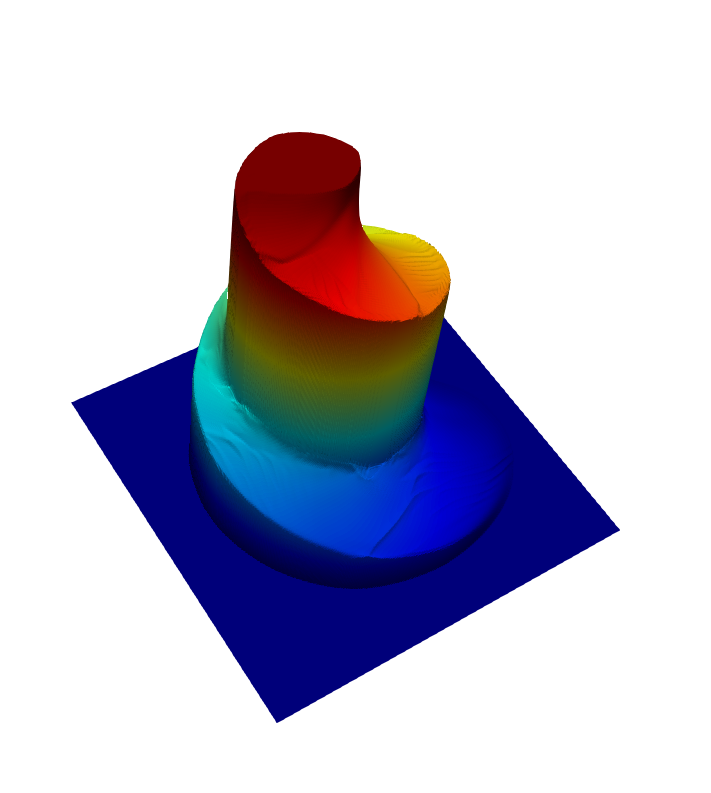}
        \caption{BP-ES}
        \label{fig:Example4-MCLES-3d}
    \end{subfigure}
    \caption{Example 4: Numerical solutions produced by different schemes at $T = 1.0$.}
    \label{fig:Example4-3d-profiles}
\end{figure}

Since the LO scheme is entropy stable, the corresponding numerical solution provides an inaccurate but qualitatively correct approximation to the spiral-shaped shock, as shown in Figure~\ref{fig:Example4-LO}. No entropy stability conditions are enforced in the BP method. As a result, we observe the merging of two shocks in Figure \ref{fig:Example4-BP}. This failure to reproduce the wave structure of the entropy solution highlights the need to apply entropy fixes. Figures \ref{fig:Example4-FCTES} and \ref{fig:Example4-MCL-ES} show that solutions produced by the entropy-stable FCT-BP-ES and BP-ES schemes preserve the correct rotating wave structure, similarly to the LO scheme. For a better visual comparison of the LO, BP-ES, and FCT-BP-ES results, we plot the corresponding solution profiles along the line connecting the top left corner of $\Omega$ to the bottom right corner. The plots are displayed in Figure \ref{fig:Example4-overline}.

\begin{figure}[h!]
    \centering
    \includegraphics[width=0.5\linewidth]{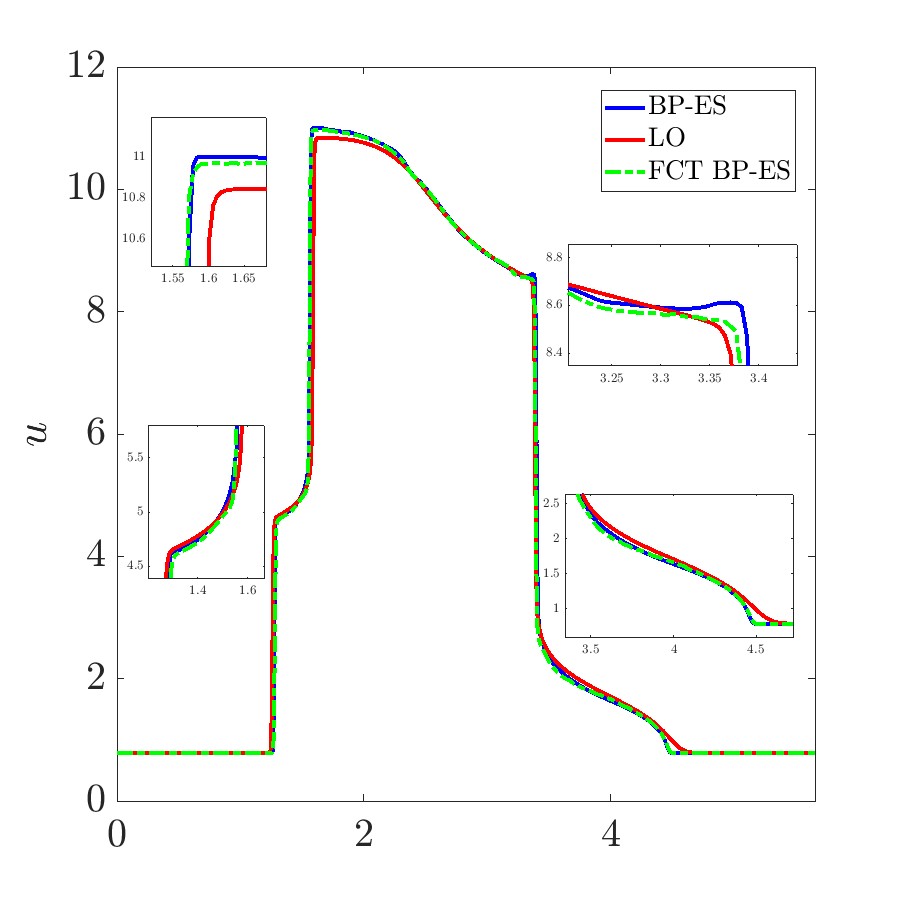}
    \caption{{Example 4. Solution profiles along the left diagonal of $\Omega$ at $T = 1.0$. }}
    \label{fig:Example4-overline}
\end{figure}

\section{Conclusions}
\label{sec:conclusions}

The local conservation property of enriched Galerkin methods makes them a promising tool for solving hyperbolic conservation laws. In the nonlinear case, preservation of local bounds and entropy stability are essential requirements for physical admissibility of finite element approximations. The proposed limiting algorithms guarantee the validity of corresponding inequality constraints and preserve optimal convergence behavior in tests with smooth exact solutions. The use of EG cell averages and CG nodal values as degrees of freedom facilitates extensions of modern limiting techniques for finite volume, DG, and CG discretizations of nonlinear hyperbolic systems, such as the Euler equations of gas dynamics and the shallow water equations \cite{kuzmin2023}. The design of customized limiters for EG discretizations of porous media flow models \cite{LEE2016,LEE201719,lee2018enriched} requires more significant effort and represents a prospective avenue for further research.

\section*{Acknowledgments}
The work of D.K. was supported by the German Research Foundation (DFG) within the framework of the priority research program SPP 2410 under grant KU 1530/30-1.
This work for S. L. and Y. Y. are supported by the National Science Foundation under Grant DMS-2208402. The authors thank Insa Schneider (TU Dortmund University) for preliminary numerical studies of the unlimited EG method in one space dimension.

\bibliographystyle{elsarticle-num} 
\bibliography{reference}

\end{document}